\documentclass[a4paper,11pt]{article}
\usepackage{cancel}
\usepackage{amsmath}
\usepackage{amssymb}
\usepackage{array}
\usepackage{ragged2e}
\usepackage{caption}
\usepackage{subcaption}
\usepackage{booktabs}
\usepackage[table,xcdraw]{xcolor}
\usepackage{url}
\usepackage{titlesec}

\usepackage{soul}
\usepackage{rotating}
\usepackage{adjustbox}
\usepackage{bm}

\usepackage{arydshln}
\usepackage{adjustbox, lipsum}
\usepackage{enumerate}
\usepackage[ruled,vlined,linesnumbered]{algorithm2e}
\usepackage{rotating}
\usepackage[utf8]{inputenc}
\usepackage{amsfonts}
\usepackage{pstricks}
\usepackage{pdflscape}
\usepackage{graphicx}
\usepackage{soul}
\usepackage{afterpage}
\usepackage{graphics}
\usepackage{cancel} 

\usepackage{rotating}
\usepackage{natbib}
\usepackage{booktabs}

\usepackage{authblk}
\usepackage{multirow}
\usepackage{mathtools}
\usepackage{scalerel}
\newtheorem{prop}{\bf Proposition}[section]

\SetKwRepeat{Do}{do}{while}

\DeclarePairedDelimiter\abs{\lvert}{\rvert}%

\newcommand{\dem}{\par \noindent{\bf Proof:} }
\newcommand{\fin}{\hfill $\square$  \par \bigskip}

\usepackage{xcolor}
\definecolor{gr}{RGB}{0,153,0}
\newcommand{\HLLP}{HLLP\xspace}
\newcommand{\ED}{ED-HLLP\xspace}
\usepackage{todonotes}
\makeatletter 
\@mparswitchfalse%
\makeatother
\normalmarginpar 
\definecolor{cyanazure}{rgb}{0, 0.7, 1}

\usepackage[normalem]{ulem}

\newcommand{\msout}[1]{\text{\sout{\ensuremath{#1}}}}
\definecolor{armygreen2}{rgb}{0.99, 0.53, 0.93}
\newcommand{\AMR}[2]{{\color{armygreen2}#1 \ifmmode\msout{#2}\else\sout{#2}\fi}}

\title{The profit-oriented hub line location problem with elastic demand}  
\date{}

\author[(a)]{Brenda Cobeña}
\author[(a)]{Ivan Contreras}
\author[(b)]{Luisa I. Mart\'inez Merino\footnote{ Corresponding author}}
\author[(b)]{Antonio M. Rodríguez-Chía}

\affil[(a)]{\small{Concordia University and Interuniversity Research Centre on Enterprise Networks,
Logistics, and Transportation (CIRRELT), Montreal, Canada, brenda.cobena@mail.concordia.ca, ivan.contreras@concordia.ca}}
\affil[(b)]{\small{Departamento de Estadística e Investigación Operativa, 
Universidad de Cádiz, 
Cádiz, Spain, luisa.martinez@uca.es, antonio.rodriguezchia@uca.es}}

\allowdisplaybreaks
\setcitestyle{authoryear,open={(},close={)}}

\makeatletter
\let\oldabs\abs
\def\abs{\@ifstar{\oldabs}{\oldabs*}}
\makeatother
\begin{document}

%
\maketitle

\begin{abstract}

 This paper deals with an extension of the hub line location problem considering demand elasticity with respect to travel times. The proposed model aims to capture the impact the hub network topology has on demand. The objective is to maximize the total revenue generated by each unit of demand using the hub line. We propose mixed-integer nonlinear formulations to model this problem. We study some properties of the nonlinear objective function associated with these formulations. Due to the inherent complexity involved in solving these nonlinear formulations with state-of-the-art solvers, we also present alternative mixed-integer linear programming formulations. Computational results compare the proposed formulations and the benefits of the presented model using benchmark instances commonly used in hub location. Moreover, a sensitivity analysis study is carried out with real data from the city of Montreal, Canada, to demonstrate the added value of incorporating demand elasticity when using the proposed model for public transportation planning. 
\end{abstract}


\section{Introduction}%
\label{sec:Introduction}%

The development of efficient and accessible public transport networks is crucial for the proper performance of mobility in metropolitan cities. Amongst others, the growth of population in urban areas \citep[][]{schafer2007long,banco2009world} and the need of eco-friendly transport modes to avoid pollution \citep{goel2017effect} are important reasons that justify why the research in urban mobility is becoming so relevant nowadays.  

In this context, hub network location models have been applied to improve the passenger mobility \citep[see,][]{Farahani,ALUMUR20211}. Generally speaking, these models address the problem of locating special facilities known as hubs and the network that establishes the connection between them. Hub facilities can be seen as transshipment nodes where there can be a change in the transport mode (e.g., bus stations and subway stations). 

In this paper, we focus on an extension of the \textit{hub line location problem} (\HLLP), a hub network design problem arising in various applications such as public transportation planning and the design of rapid transit systems. In these applications, hubs (stops, tram stations, etc.) are not fully interconnected, and direct connections between origin-destination (OD) pairs are allowed. A path structure at the hub-level network fits with these applications. We refer to  \cite{LeeRo}, \cite{LabbeYaman}, \cite{Yaman}, \cite{DeSaCamagoMiranda}, \cite{ContrerasTanash}, and \cite{Contreras2021}, for additional hub network design models considering other network topologies. 

The HLLP was first introduced by \cite{DeSa2015} and seeks to locate $p$ hubs connected by a path (or line) composed of a set of $p-1$ hub edges. It minimizes the sum of the total weighted travel time between OD pairs. The flows represent passengers or users travelling between OD pairs who wish to minimize their commute time. They will use the hub line whenever time savings are perceived. Otherwise, passengers will use a direct link. The HLLP incorporate other aspects relevant to model travel times such as the access and exit times incurred when using the hub network. 

The HLLP assumes that demand is inelastic, i.e., it does not depend on the design of the resulting hub network. However, this assumption may be unrealistic.  For instance, in the case of subway networks, when a new line is opened, an increase in the demand flow between newly connected neighborhoods can be perceived since the reduction of travel time becomes attractive for some users. Therefore, if demand elasticity is not considered during network planning, the long-term benefits of the hub line hub may not be correctly captured in the model. 

In the literature, we can find some models where demand elasticity is one of the factors to take into consideration.  \cite{DreznerDrezner2001}  was the first paper applying the gravity rule to hub location models. On the other hand, in competitive facility location, some problems maximize market capture considering elastic demand \citep[see,][]{Marianov99}. 
\cite{Marianov09} deals with a competitive hub location problem where customers have gravity-like utility functions. \cite{Marianov2005} and \cite{Marianov2008} propose a model for facility location where elasticity is modeled by a supply-demand equilibrium equation. Demand elasticity is considered in other location problems. For instance, the profit-maximizing service network design problem introduced in \cite{Aboolian2012}, the model presented in \cite{Kuiteing2017} related to network pricing, and profit-oriented multi-commodity network design with elastic demand \citep{Zetina2019}.

In this paper we extend the HLLP by using a gravity model to incorporate demand elasticity with respect to travel times with a profit maximization perspective. This extension, denoted as the \textit{profit-oriented hub line location problem with elastic demand} (\ED), is more realistic since, in general, passenger demands change according to population mobility after establishing the line. Therefore, the \ED can be used to provide new insights to public stakeholders such as network designers, city planners, and public transport managers. The considered gravity model is based on Newton's fundamental law of attraction. It measures interactions between all the possible location pairs by multiplying their attributes over the distances between them. Moreover, gravity models have been used as the basis to model many application areas such as understanding the effect of migration on voting patterns, planning a new transportation service (or a new road), determining the optimal size of a shopping development, defining retail shopping boundaries, and traffic flows \citep[see,][]{Haynes}. 

In addition to introducing the \ED, the main contributions of this work are the following. We present two main mixed-integer nonlinear programming formulations (MINLP) for the \ED that use arc-based variables to model OD paths. The main difference between these formulations is the way in which the hub line topology is modeled. Given that these formulations use a nonlinear objective, which is difficult to optimize with state-of-the-art MINLP solvers, we also propose three mixed-integer linear programming formulations (MILP) that use path-based variables to model OD paths. We exploit the considered network topology to develop an efficient algorithm to generate all candidate OD paths. This allow us to solve the resulting compact MILP formulations directly with commercial solvers for small to medium size instances without the need of using more sophisticated branch-and-price algorithms. To evaluate the efficiency and limitations of all MINLP and MILP formulations, computational experiments were performed on benchmark instances with up to 25 nodes and $p=7$. Moreover, a sensitivity analysis was carried out with real data from the city of Montreal, Canada, to demonstrate the added value of incorporating demand elasticity in public transportation planning. 

The remainder of the paper is organized as follows. Section \ref{notation} provides a formal definition of the problem. Section \ref{nonlinear} describes the MINLP formulations for the \ED and analyzes some properties of the nonlinear objective function considered in these formulations. Section \ref{linear} presents alternative MILP formulations for the \ED and the algorithm developed to generate all candidate OD paths. The results of computational experiments are reported in Section \ref{computational}. Conclusions follow in Section \ref{conclusions}.

\section{Problem definition and notation}\label{notation}

Let $G=(N,E)$ be a graph, where $N$ is the set of nodes and $E$ the set of edges $e:=[k,m]$ with $k<m$, $k,m\in N$. We also consider $\widetilde{G}=(N,A)$ which is the directed version of graph $G$, i.e., $A=\{(k,m)\cup (m,k):e=[k,m]\in E\}$ is the set of arcs induced by $E$. 
Besides, let $C$ be a set of OD pairs whose demand must be routed either using a hub line or directly from origin to destination. Each OD pair will be referred to as commodity $c\in C$, its origin is denoted by $o_c$ and its destination as $d_c$.
{The travel time of arc $(i,j)\in A$, represented by $t_{ij} >0$, is defined as the
shortest time required to travel from $i$ to $j$, without using the established hub line, i.e., using the original (physical) network. Note that this
definition of travel times leads to $t_{ij}$ values that satisfy the triangle inequality property.}
 For each $i\in N$, $P_i$ is a weight related to its population size. 

When using the hub line, the access time to enter the first hub $k$ is given by $\tilde{t}_k^a$. Similarly, the exit time from the last hub $m$ is denoted by $\tilde{t}_m^e$.  The access time $\tilde{t}_k^a$ incorporates both the time required to change the mode of transportation between an access arc and a hub arc at node $k$ and also the average waiting time to access hub $k$. The exit time $\tilde{t}_m^e$ includes the time required to change the mode of transportation at the last hub node $m$. 
{Because of the triangle inequality property of travel
times $t_{ij}$, there exists a solution network of the HLLP
that routes demands $w_{ij}$ either with a direct connection
between $i$ and $j$ or with a path containing at most
two access arcs and at least two hub nodes and one
hub arc.} Thus, once a commodity leaves the hub line, it cannot access the hub line again.

The \ED  seeks to maximize the time savings related to the total flow of the set of commodities. 
Each unit of time reduction for a commodity $c\in C$ due to the use of the hub line has an associated revenue $R_c\geq 0$. In addition, the total time savings depend on the final demand associated with the commodities routed by the hub line. Consequently, the modeling of the elastic demand has an important impact on the considered problem.

Gravity models consider the attraction between two locations as directly proportional to their size and inversely proportional to the distance between them. The gravity distribution model has the following form \citep[][]{taaffe1996geography}:
\begin{align}\label{gravity1}
W_{i j} = \frac{P_i P_j}{(d_{i j})^{r}},
\end{align}
where,  $W_{i j}$ is the demand between $i$ and $j$. ${P_i}$ and ${P_j}$ are features related to locations $i$ and $j$, respectively. Usually, these values are obtained by considering the population, number of available jobs or the gross domestic product of the location. $d_{i j}$ is the distance or travel time between the origin $i$ and the destination $j$. Finally, ${r}$ is the exponent of distance (time), and it is associated with the efficiency of the transport system between two locations. In \eqref{gravity1}, the initial population  ${P_i}$  is usually referred to as $production$ whereas the second population ${P_j}$ as $attraction$ \citep[see,][]{de2011modelling}. 

In the context of the \ED, the above gravity model determines the final demand of a commodity $c$, denoted as $w_c$, as follows:
$$w_c=\frac{P_{o_c}P_{d_c}}{(T_{c})^r}.$$

The product in the numerator can be seen as the initial attractiveness between $o_c$ and $d_c$. $T_c$  denotes the final travel time between $o_c$ and $d_c$, either using the hub line or without using it. Once the hub line is located, if a commodity does not save time using the hub line, its associated demand will be 
$$w_c=\frac{P_{o_c}P_{d_c}}{(t_{o_c d_c})^r}.$$
Otherwise, the final demand of $c$ will depend on the resulting shortest time $(t'_{o_c d_c})$ between $o_c$ and $d_c$ using the hub line:
$$w_c=\frac{P_{o_c}P_{d_c}}{(t'_{o_c d_c})^r}.$$

Different values of $r$ can be considered depending on the particular application of the gravity model. Usually, $r$ ranges between 0.5 and 2 \citep[see,][]{FKelly1989}. Table \ref{tab1}  summarizes the notation introduced for this problem.

\begin{table}[htbp]
	\centering
	\begin{tabular}{ll}
		\toprule
		Notation & Description \\
		\midrule
		$G$    &  Graph where $N$ is the set of nodes and $E$ the set of edges.  \\
		$\widetilde{G}$&  Directed graph induced by $G$.\\
		$C$     & Set of commodities. \\
		$o_c$ & Origin of commodity $c\in C$.\\
		$d_c$ & Destination of commodity $c\in C$.\\
		$t_{o_c d_c}$ & Shortest time required to travel from $o_c$ to $d_c$  in $G$ (without using \\
		& the hub line.)\\
		$\tilde{t}^a_i$ & Access time to the hub line through node $i\in N$.\\
		$\tilde{t}^e_i$ & Exit time from the hub line through node $i\in N$.\\
		 $R_c$ & Revenue for each unit of time reduction associated with $c\in C$.\\
		 $P_i$& Weight related to population size of $i\in N$. \\
		 $r$ &Value of the exponent considered in the gravity model.\\ 
		\bottomrule
	\end{tabular}%
    \caption{Overview of \ED parameters. }
	\label{tab1}%
\end{table}%

\section{MINLP formulations for the \ED}
\label{nonlinear}

In this section, we first introduce three MINLP formulations for the \ED inspired on the \HLLP formulation presented in \cite{DeSa2015}. We then present, in Section \ref{objectiveprop}, some properties of their objective function. Finally, in Section \ref{impro_nolinear} we provide some preprocessing, valid inequalities, and fixing-variables criteria in order to improve the performance time of these formulations. 

In the first formulation, denoted as $F1$, we consider the following decision variables:

\begin{eqnarray*}
	z_k &=& \begin{cases}
		1, & \mbox{if  a hub is located at node $k$,} \\
		0, & \mbox{otherwise,}
	\end{cases}
	\quad k \in N, \\
	y_{e} &= & \begin{cases}
		1,& \mbox{if a hub edge is located in edge $e$ of graph $G$}, \\
		0, & \mbox{otherwise},
	\end{cases}
\quad e\in E.\\
\end{eqnarray*}
In $F1$, $\bm{z}$- and $\bm{y}$-variables determine the nodes and the edges that will define the hub line topology. We also need to define the following sets of decision variables to model OD paths:
\begin{eqnarray*}
	a_{ck}&= &\mbox{fraction of commodity $c$  that enters the hub line through hub $k$},\\
	&&c\in C,k \in N,\\
	b_{cm}&= &\mbox{fraction of commodity $c$  that exits the hub line through hub $m$},\\
	&&c\in C, m \in N,\\
	x_{ckm} &= & \mbox{fraction of commodity $c$ routed from hub node $k$ to hub node $m$,}\\
	&& \mbox{for $c\in C,\, k,m\in N$,}\\
	e_{c} &=& \begin{cases}
	1,&\mbox{if commodity $c$ is sent directly from $o_c$ to $d_c$,}\\
	0,&\mbox{otherwise,}\\
	\end{cases}c\in C.
\end{eqnarray*}

We also need to use a set of variables to model the connectivity of the hub line:
\begin{equation}
    f_{km}=\mbox{flow from node $k$ to node $m$ with $k,m\in N$.}
    \label{flowvar}
\end{equation}

Finally, we include a set of variables to represent the overall travel time associated with each commodity (either using the hub line or traveling directly from origin to destination):
\begin{eqnarray*}
t'_c&=& \mbox{total travel time between $o_c$ and $d_c$,}\,c\in C.
\end{eqnarray*}
Using the aforementioned variables, the \ED can be stated as:

\begin{alignat}{3}
(F1)\quad& \max&  \quad & \sum_{c\in C}
R_c\frac{P_{o_c}P_{d_c}}{(t'_c)^r}\left(t_{o_c d_c}-t'_c\right)
\nonumber\\
&\mbox{s.t.}&&t'_c=\sum_{k\in N}(t_{o_c k}+\tilde{t}_k^a)a_{ck}+\sum_{k\in N}\sum_{\stackrel{m\in N}{ (k,m)\in A}}\alpha_{km}t_{km}x_{ckm}+\nonumber\\
&&&\hspace*{0.7cm}\sum_{m\in N}(t_{ m d_c}+\tilde{t}_m^e)b_{cm}+t_{o_c d_c}e_c &\quad&\hspace{-1.6cm}c\in C,\label{finalt}\\
& &&\sum_{k\in N}a_{ck}+e_c=1 &&\hspace{-1.6cm}c\in C,\label{a_rest}\\
& &&\sum_{m\in N}b_{cm}+e_c=1 &&\hspace{-1.6cm}c\in C,\label{b_rest}\\
&&&a_{ck}+\sum_{\stackrel{m\in N}{m\neq k}}x_{cmk}=b_{ck}+\sum_{\stackrel{m\in N}{m\neq k}} x_{ckm} &\quad&\hspace{-1.6cm}c\in C,k\in N,\label{ab_rest}\\
&&&a_{ck}\leq z_k &\quad&\hspace{-1.6cm}k\in N,c\in C,\label{az_rest}\\
&&&b_{cm}\leq z_m &\quad&\hspace{-1.6cm}m\in N,c\in C,\label{bz_rest}\\
&&&x_{ckm}+x_{cmk}\leq y_e &\quad&\hspace{-1.6cm}c\in C,e=[k,m]\in E,\label{xy_rest}\\
&&&\sum_{i\in N}z_k=p &\quad&\label{p_rest}\\
&&&\sum_{e\in E}y_e=p-1 &&\label{yp_rest}\\
&&&\sum_{\stackrel{m\in N}{[k,m]\in E}}y_{[k,m]}+\sum_{\stackrel{m\in N}{[m,k]\in E}}y_{[m,k]}\leq 2z_{k} &\quad\hspace{-1.5cm}&\hspace{-1.6cm}k\in N,\label{yz_rest}\\
&&&\sum_{\stackrel{m\in N}{(k,m)\in A}}f_{km}\leq (p-1)z_k &&\hspace{-1.6cm} k\in N,\label{flow1}\\
&&&\sum_{\stackrel{j\in N}{(j,m)\in A}}f_{jm}-\sum_{\stackrel{j\in N}{(m,j)\in A}}f_{mj}\geq z_m+(z_k-1)p &&\hspace{-1.5cm}
m,k\in N, m<k,\label{flow2}\\
&&&f_{km}+f_{mk}\leq (p-1)y_{e} &&  \hspace{-1.6cm} e=[k,m]\in E\label{flow3}\\
&&&y_{e}\in\{0,1\} &\quad&\hspace{-1.6cm}e\in E,\label{y_dom}\\
&&&z_k\in\{0,1\} &\quad&\hspace{-1.6cm}k\in N,\label{z_dom}\\
&&&e_c\in \{0,1\} &\quad&\hspace{-1.6cm}c\in C,\label{e_dom}\\
&&&x_{ckm}\geq 0 &\quad&\hspace{-1.6cm}c\in C,(k,m)\in A,\\
&&&f_{km}\geq 0 &\quad&\hspace{-1.6cm}(k,m)\in A,\label{f_dom}\\
&&&a_{ck}\geq 0 &\quad&\hspace{-1.6cm}c\in C,k\in N,\\
&&&b_{cm}\geq 0 &\quad&\hspace{-1.6cm}c\in C,m\in N.\label{b_dom}
\end{alignat}

The objective function maximizes the total revenue in terms of travel time reduction for delivering the realized demand. Constraints \eqref{finalt} compute the travel time for each commodity, either through the hub line or without using it. Constraints \eqref{a_rest} and \eqref{b_rest} ensure that commodity $c$ is either delivered through a hub line or directly.
Constraints \eqref{ab_rest} are the flow conservation constraints of the commodities at the hubs of the hub line. Constraints \eqref{az_rest} and \eqref{bz_rest} restrict that only open hubs can be the access and exit hub nodes for a commodity, respectively. In the hub line, a commodity can only be transported using an open hub arc, this is specified in constraints \eqref{xy_rest}. Constraints \eqref{p_rest} and \eqref{yp_rest} ensure that $p$ hubs and $p-1$ hub edges are opened, respectively. Constraints \eqref{yz_rest} ensure that the edges of the hub lines have hubs as end-nodes. Moreover, \eqref{flow1}--\eqref{flow3} are flow constraints that avoid subtours and ensure that the resulting hub structure is connected. 

An alternative formulation for \ED can be obtained by modeling connectivity and subtour elimination constraints through a path over $\widetilde{G}$, in such a way that the resulting hub line will be this path, replacing arcs by edges. Hence, we define the following variables:

\begin{eqnarray}
	y'_{km}&=&\begin{cases}
		1,&\mbox{ if arc $(k,m)$ belongs to the hub line,}\\
		0,&\mbox{otherwise,}\\
	\end{cases}
	\label{y'var}
\end{eqnarray}
for each $(k,m)\in A$. It is necessary to include the following set of auxiliary variables: $l_i\geq 0$, for $k\in N$. Observe that the idea behind using $\bm{y'}$-variables is to take advantage of
using a path in the directed graph $\widetilde{G}$ as hub line structure, although the commodities can be supplied in both directions, i.e., if $[i,j]\in E$ such that $y'_{ij}=1$ or $y'_{ji}=1$, then edge $[i,j]$ is a hub edge. Using these new sets of variables, the \ED can be stated as:

\begin{alignat}{3}
(F2)&\quad \max&  \quad &\sum_{c\in C}
R_c\frac{P_{o_c}P_{d_c}}{(t'_c)^r}\left(t_{o_c d_c}-t'_c\right)\nonumber\\
&\quad\mbox{s.t.}&&\eqref{finalt}-\eqref{bz_rest},\eqref{p_rest},\eqref{z_dom}-\eqref{b_dom},\nonumber\\
&&&x_{ckm}+x_{cmk}\leq y'_{km}+y'_{mk} &\quad&c\in C,[k,m]\in E,\label{xy_rest2}\\
&&&\sum_{k\in N}\sum_{\stackrel{m\in N}{(k,m)\in A}}y'_{km}=p-1,&&\label{yp_rest2}\\
& &&  \sum_{\stackrel{m\in N}{(k,m)\in A} } y'_{km} \le  z_k &\quad&k \in  N, \label{yz_1}\\
& && \sum_{\stackrel{m \in N}{(m,k)\in A}} y'_{mk}\le z_k &\quad&  k \in  N, \label{yz_2}\\  
&&&l_k-l_m+n y'_{km}\leq n-1 &\quad&(k,m)\in A,
\label{MTZ}\\
&&&y'_{km}\in\{0,1\} &\quad&(k,m)\in A.\label{y_dom2}
\end{alignat}

Similar to \eqref{xy_rest}, constraints \eqref{xy_rest2} ensure that only open hub edges can be used to route a commodity. Constraints \eqref{yz_1} and \eqref{yz_2} ensure that there is at most one incoming arc and one outgoing arc of the path in each hub node. In $F2$, we need a family of constraints to avoid subtours and ensure the connectivity of the hub line. The so-called Miller-Tucker-Zemlin(MTZ) constraints have been traditionally used to guarantee connectivity and to avoid subtours in routing models \citep[see,][]{MTZ,Laporte,Gouveia1995,Bektas,LandeteMarin}. In 
$F2$, constraints \eqref{MTZ} are MTZ constraints.

Note that constraints \eqref{yz_1} and \eqref{yz_2} can be aggregated resulting, in the following formulation:
\begin{alignat}{3}
(F2')&\quad \max&  \quad &\sum_{c\in C}
R_c\frac{P_{o_c}P_{d_c}}{(t'_c)^r}\left(t_{o_c d_c}-t'_c\right)\nonumber\\
&\quad\mbox{s.t.}&&\eqref{finalt}-\eqref{bz_rest},\eqref{p_rest},\eqref{z_dom}-\eqref{b_dom},\eqref{xy_rest2},\eqref{yp_rest2},\eqref{MTZ},\eqref{y_dom2},\nonumber\\
& &&  \sum_{\stackrel{m\in N}{(k,m)\in A}} y'_{km}+ \sum_{\stackrel{m \in N}{(m,k)\in A}} y'_{mk}\le 2 z_k,&\quad&k \in  N. \label{yz_new}
\end{alignat}

Moreover, the following set of valid inequalities can be added to $F2'$:
\begin{eqnarray}
\sum_{\stackrel{l\in N}{ (k,l)\in A}} y'_{kl}\geq z_k+z_m-1,&& k\in N\setminus \{n\},m=k+1,\ldots,n.\label{new}
\end{eqnarray}

Constraints \eqref{new} in $F2'$ ensure that there is at least one outgoing hub arc from each hub node of the path, except the one with the largest index. The idea of using constraints \eqref{MTZ}, \eqref{yz_new} and \eqref{new} is to obtain an oriented path such that the hub with the biggest index does not have an outgoing hub arc. This avoids symmetries in feasible solutions. For instance, a possible solution could be $y'_{11,7}=y'_{7,16}=y'_{13,16}=1$. As a result, the edges of the hub line are $[7,11]$, $[7,16]$ and $[13,16]$ and they can be traversed in both directions. Observe that valid inequalities \eqref{new} cannot be applied to formulation $(F2)$   since constraints \eqref{yz_1} and \eqref{yz_2} avoid some feasible solutions. Following the previous mentioned example, using constraints \eqref{yz_1} and \eqref{yz_2}  would avoid the case $y'_{7,16}=y'_{13,16}=1$.

\subsection{Properties of the Objective Function} \label{objectiveprop}

The objective functions in formulation $(F1)$, $(F2)$ and $(F2')$ are the same, and we can derive the following properties of this nonlinear function.

\begin{prop}\label{prop1}
	The objective function of $(F1)$, $(F2)$ and $(F2')$ is nonincreasing  and convex on $\bm{t'}=(t'_1,t'_2,\ldots,t'_{|c|})$.
\end{prop}
\dem
Observe that the mentioned objective function is given by,

$$f({\bm{t'}})=\sum_{c\in C}
R_c\frac{P_{o_c}P_{d_c}}{(t'_c)^r}\left(t_{o_c d_c}-t'_c\right),$$

where $$t'_c=\sum_{k\in N}(t_{o_c k}+\tilde{t}_k^a)a_{ck}+\sum_{k\in N}\sum_{\stackrel{m\in N}{ (k,m)\in A}}\alpha_{km}t_{km}x_{ckm}+\sum_{m\in N}(t_{m d_c}+\tilde{t}_m^e)b_{cm}+t_{o_c d_c}e_c.$$
Taking into account that the addends of the objective function with $e_c=1$ will take the null value due to 
constraints \eqref{a_rest}-\eqref{ab_rest}, then 

$$f({\bm{t'}})=\sum_{\stackrel{c\in C}{e_c=0}}
R_c\frac{P_{o_c}P_{d_c}}{(t'_c)^r}\left(t_{o_c d_c}-t'_c\right),$$

We now analyze each addend of the function above. Let $f_c({\bm t'_c})=R_c\frac{P_{o_c}P_{d_c}}{(t'_c)^r}\left(t_{o_c d_c}-t'_c\right)$. We obtain that

\begin{itemize}
	\item $f_c'( t'_c)= R_c P_{o_c}P_{d_c}\frac{ (r-1)t'_c- r t_{o_c d_c}}{{t'}_c^{r+1}}$. Observe that
	$f'_c(t'_c) \leq 0$,  since $t'_c\leq t_{o_c d_c}$, and consequently, $(r-1) t'_c \le r t_{o_c d_c}$.
	\item $f_c''(t'_c)= R_c P_{o_c}P_{d_c} r \frac{ (r+1) t_{o_c d_c} - (r-1) t'_c}{{t'}_c^{r+2}}$. Observe that
	$f''(t'_c)\geq 0$. 
\end{itemize}
Then, the objective function is nonincreasing and convex.
\fin

\subsection{Preprocessing, variable fixing and valid inequalities}\label{impro_nolinear}

In the previous section, we proposed several formulations for the \ED: $(F1)$, $(F2)$, $(F2')$ and $(F2')$+\eqref{new}. All of them are nonlinear due to their objective function. Consequently, in Section \ref{computational}, a general-purpose global optimization solver has been used to solve these models. In particular, we used BARON 21.1.13 solver \citep[][]{Kilin2018ExploitingII,Nohra} through AMPL modeling language. In this global optimization solver, it is necessary to include adequate upper and lower bounds for each nonlinear term appearing in the objective function. In order to provide good lower and upper bounds on $f_c( t'_c)=R_c\frac{P_{o_c}P_{d_c}}{(t'_c)^r}\left(t_{o_c d_c}-t'_c\right)$, where $t'_c$ is the shortest travel time to transport commodity $c$ using the resulting hub line, we provide the following procedure.  

Note that $0$ is a valid lower bound for $f_c(t'_c)$. Besides, from Proposition \ref{prop1}, $f_c(t'_c)$ is nonincreasing in $t'_c$. Consequently, if the smallest value of $t'_c$ is $t_0$, then we have that $f_c(t_0)$ is an upper bound for $t'_c$. Therefore, in order to obtain an upper bound for this nonlinear term, we look for the smallest value of $t'_c$, i.e., the shortest path between $o_c$ and $d_c$ either using the hub line or not.

For each $c\in C$, we create an {auxiliary} directed graph $G_c$ that allows us to build all possible paths between $o_c$ and $d_c$ using the hub line with an associated travel time smaller than or equal to $t_{o_c d_c}$. This graph $G_c$ is composed by a set of nodes $N_c=N\cup\{o'_c,d'_c\}$, i.e., we have the original set of nodes, one copy of the origin related to commodity $c$, $o'_c$, and another copy of destination associated with commodity $c$, $d'_c$. Consequently, $t_{o_c o'_c}=0$ and $t_{d'_c d_c}=0$. The arcs associated with $G_c$, denoted as $A_c$, and  the travel time of using each arc (denoted as $t_{ij}^c$ for $(i,j)\in A_c$) can be obtained by Algorithm \ref{alg-1}. This procedure ensures that we are only considering the arcs that allow us to obtain paths with a travel time smaller than or equal to the travel time in the original graph.

\begin{algorithm}[H]
	\label{alg-1}
	\KwData{Graph $G=(N,E)$ with $t_{ij}$ for $e=[i,j]\in E$ and $c\in C$.} 
	\KwResult{{Auxiliary} graph $G_c=(N_c,A_c)$ and $t^c$.} 
		$N_c=N\cup \{o'_c,d'_c\}$.\\ 
		$A_c$ and $t^c$ are generated as follows:\\
		\For{$i,j\in N_c$}{
			\If{$i\neq j$ and $j\neq o_c$ and $i\neq d_c$ and $i\neq d'_c$ and $j\neq o'_c$}{
			    \If{$t_{o_ci} + \tilde{t}^a_{i} + \alpha\cdot t_{ij}  \le t_{o_c d_c} $}{
			    include $(o_c, i)$ in $A_c$,\\
				$t^c_{o_ci}= t_{o_ci} + \tilde{t}^a_{i}$}
				\If{$  \alpha\cdot t_{ij} +  t_{j d_c} + \tilde{t}^e_{j}  \le t_{o_c d_c} $}{
				include $(j,d_c)$ in $A_c$,\\
				$t^c_{j d_c}= t_{j d_c} + \tilde{t}^e_{j}$}
				\If{$  \alpha\cdot t_{ij} \le t_{o_c d_c} $}{
				include $(i,j)$ in $A_c$,\\
				$t^c_{ij}= \alpha\cdot t_{ij}$}
}
}
	\For{$i\in N_c\setminus\{o_c,d_c,o'_c,d'_c\}$}{
			\If{degree of $i$ is smaller than or equal to 1 }{
			$N_c=N_c\setminus\{i\}$
	}}
\Return $G_c=(N_c,A_c)$ and $t^c$.
	\caption{Obtaining {auxiliary} graph $G_c=(N_c,A_c)$ and travel times $t^c$ of each arc for  $c\in C$.}
\end{algorithm}

Note that, after obtaining $A_c$, the set of considered nodes can also be reduced. If a node $i\in N\setminus\{o_c,d_c\}$ has a degree smaller than or equal to $1$, it means that there is not path between $o_c$ and $d_c$ with a related time smaller than or equal to $t_{o_c d_c}$. Consequently, if degree of $i$ is smaller than or equal to $1$, then $i\notin N_c$. This is described in lines 14-16 of Algorithm \ref{alg-1}.

\begin{figure}[!htb]
	\centering
	\begin{subfigure}[b]{0.4\textwidth}
		\centering
		\includegraphics[width=0.6\linewidth]{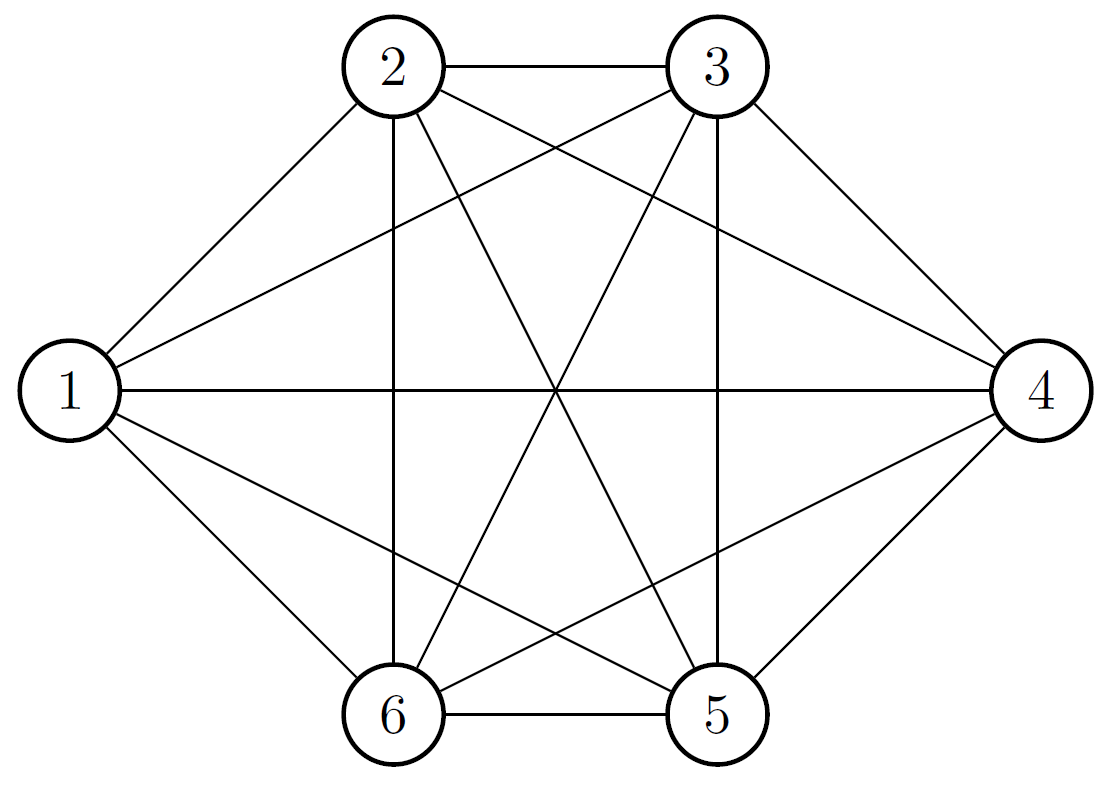}
		\caption{Example of an original\\ graph $G$}
		\label{ejemplo2}
	\end{subfigure}%
	\begin{subfigure}[b]{0.4\textwidth}
		\centering
		\includegraphics[width=1\linewidth]{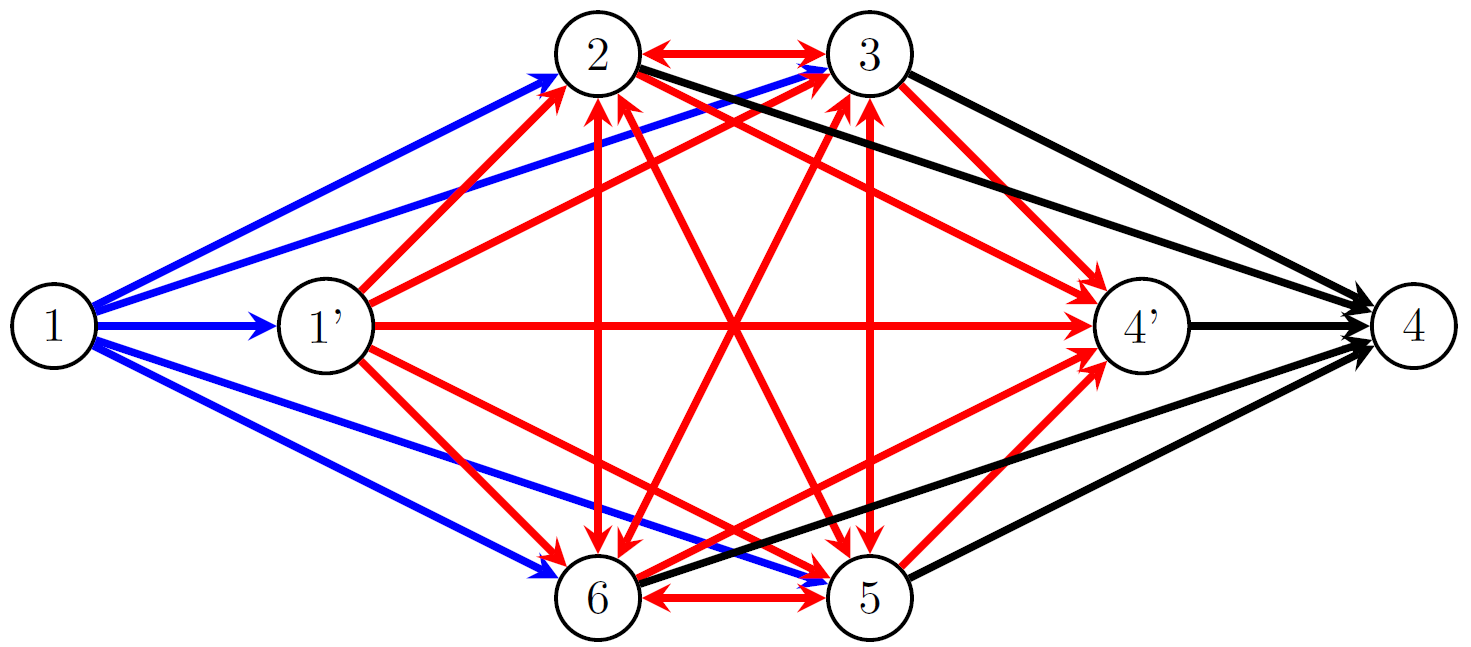}
		\caption{{auxiliary} network $G_c$ of $G$ to obtain the hub line paths from $o_c=1$ to $d_c=4$}
		\label{ejemplo3}
	\end{subfigure}
	\caption{Construction of auxiliary graph $G_c$ for $o_c=1$ and $d_c=4$}
\end{figure}

Observe that for those paths that use the arc $(o_c,o'_c)$, the access node to the hub line is $o_c$. Similarly, if arc $(d'_c,d_c)$ belongs to the path, then the path uses node $d_c$ as exit node of the hub line. An example of using an auxiliary graph $G_c$ is illustrated in Figures \ref{ejemplo2} and \ref{ejemplo3}. Note that, in this example, we are considering a complete graph of six nodes that is presented in Figure \ref{ejemplo2}. Figure \ref{ejemplo3} shows the transformation of original edges in arcs to obtain an auxiliary graph that allows to obtain all possible hub line paths from node 1 to node 4. Observe that, in this example, we assume that $t^c_{ij}\leq t_{o_c d_c}$ for $i,j\in N_c$.

To obtain the shortest time path of commodity $c$, we use the function  {\sl shortest\_simple\_paths} from Python package {\sl Networkx}. If the obtained path has hub nodes smaller than or equal to $p$, then the resulting path is the shortest path of commodity $c$ using a hub line. Therefore, if the time associated with the path is denoted as $t^{LB}_c$, we have that 
\begin{equation}
0\leq f_c(t'_c)\leq f_c(t^{LB}_c),\quad \mbox{ for }c\in C.\label{lb_ub}
\end{equation}

If the path has a number of hub nodes greater than $p$, then we continue looking for the path with the second-shortest time. In Algorithm \ref{alg0}, we specify the details of the implementation. Constraints \eqref{lb_ub} must be included in $F1$, $F2$, $F2'$ and $F2'+$\eqref{new} in order to be solved by BARON.

	\begin{algorithm}[H]
		\label{alg0}
		\KwData{Graph $G=(N,E)$ with $t_{ij}$, for $i,j\in N$.
		A commodity $c\in C$. }
		\KwResult{An upper bound for $f_c(t'_c)$.}
		Create a $G_c=(N_c,A_c)$ with travel times $t^c$ following Algorithm \ref{alg-1}.\\
		\For{path in shortest\_simple\_paths($G_c$,$o_c$,$d_c$,$t^c$)}{
		\If{time(path)$< t_{o_c d_c}$}{$f_t(t^{LB}_c)=R_c \frac{P_{o_c}P_{d_c}}{\mbox{time(path)}^r}\left(t_{o_cd_c}-\mbox{time(path)}\right).$\\
		\If{(number of nodes of path $-2)\leq p$}{break}
		}
		\Else{$f_c(t^{LB}_c)=0$; break;\\}
		}
		\Return{ $f_c(t^{LB}_c)$ is an upper bound for $f_c(t'_c)$.}
		\caption{Given $c\in C$, obtain an upper bound for $f_c(t'_c)$.}
	\end{algorithm}

Moreover, in $F1$, $F2$, $F2'$ and $F2'+$\eqref{new} we can fix some variables and add some valid inequalities as follows:

\begin{itemize}
	\item[i)] For $c \in C$, $k\in N$, if $t_{o_c k}+\tilde{t}^a_k>t_{o_cd_c}$, then any path with $(o_c,k)$ imply a travel time greater than the time of direct path O/D. Consequently, 
	\begin{eqnarray}
	a_{ck}&=&0.\label{fix0}
	\end{eqnarray}
	\item[ii)] Similarly, for $c \in C$, $k\in N$, if $t_{d_c k}+\tilde{t}^e_k>t_{o_cd_c}$, then
		\begin{eqnarray}
	b_{ck}&=&0.\label{fix1}
	\end{eqnarray}
	\item[iii)]For $c\in C$, the access node for a path in the hub line cannot be the destination node:
	\begin{eqnarray}
	a_{cd_c}=0.
	\end{eqnarray}
	\item[iv)] Similarly, for $c\in C$, the exit node for a path in the hub line cannot be the origin node:
	\begin{eqnarray}
	b_{c o_c}=0.
	\end{eqnarray}
	\item[v)] For each $c\in C$, $[k,m]\in E$, if $\alpha t_{km}>t_{o_c d_c}$, then there is no path using edge $[k,m]$ that implies a time for commodity $c$ smaller than the direct path. Consequently,
	\begin{eqnarray}
	x_{ckm}=0.
	\end{eqnarray}
	\item[vi)] For each $c\in C$, $k,m\in N$, $k\neq m$, if $t_{o_c k}+\tilde{t}^a_k+\tilde{t}^e_m+t_{m d_c}>t_{o_c d_c}$, then we can add the following valid inequality,
	\begin{eqnarray}
	a_{ck}+b_{cm}\leq 1.\label{fixn}
	\end{eqnarray}
\end{itemize}

\section{MILP formulations for the \ED}
\label{linear}

 In this section we present MILP formulations for the \ED that exploit the fact, if the travel times are known a priori, the objective function becomes linear. Given a commodity $c\in C$, let $\mathcal{P}_c$ denote the set of all possible paths  in the original graph using a hub line of $p$ hubs with an associated travel time smaller than $t_{o_c d_c}$. The reader is referred to Subsection  \ref{paths} where an efficient procedure to obtain all these paths for each commodity is developed. 
 
 Each path $\pi\in\mathcal{P}_c$ can be expressed as:
$$\pi=[o_c,h_1,\ldots,h_k,d_c],$$
where $h_m$, for $m=1,\ldots k$ with $k\leq p$, denote the hubs that path $\pi$ traverses in its correct order. Specifically, $h_1$ and $h_k$ are the access and exit nodes of the hub line, respectively. Note that for some paths where origin or destination are hubs, we have that $h_1=o_c$ or $h_k=d_c$. Then, the travel time for routing commodity $c\in C$ via path $\pi\in \mathcal{P}_c$ ($\tau_{\pi c}$) is: 
\begin{equation*}
\tau_{\pi c}=t_{o_c h_1}+\tilde{t}^a_{h_1}+\sum_{m=1}^{k-1}\alpha t_{h_m h_{m+1}}+\tilde{t}^e_{h_k}+t_{h_k d_c}.
\end{equation*}

Therefore, for each $\pi\in\mathcal{P}_c$, its associated travel time $\tau_{\pi c}$ is known and, consequently, the profit of path $\pi$ for commodity $c$ using the gravity model is:
\begin{equation*}
    g'_{\pi c}= R_c\frac{P_{o_c} P_{d_c}}{(\tau_{\pi c})^r}(t_{o_c d_c}-\tau_{\pi c}).
\end{equation*}

Since this profit can be calculated for each $\pi\in \mathcal{P}_c$ for $c\in C$, a linear objective function can be used to model the \ED. 
In order to do this, we need to define the following family of variables:
\begin{eqnarray*}
	v_{\pi c}&=&\begin{cases}
		1, & \mbox{if commodity $c$  is delivered using  path $\pi$},\\
		0, & \mbox{otherwise,}
	\end{cases}\quad \pi \in \mathcal{P}_c,c\in C,
\end{eqnarray*}

and the following parameters:
$$h_{e}^{\pi c}=\begin{cases}
1,&\mbox{if path $\pi\in\mathcal{P}_c$ contains the arc $(k,m)$ or $(m,k)$ defined by }\\
&\mbox{edge $e=[k,m]\in E$,}\\
0,&\mbox{otherwise.}\\
\end{cases}$$

 The \ED can be stated as the following MILP:
\begin{alignat}{3}
(F1_L)\quad& \max&  \quad & \sum_{c\in C}\sum_{\pi\in \mathcal{P}_c} g'_{\pi c}v_{\pi c}\nonumber\\
&\mbox{s.t.}&&\eqref{p_rest}-\eqref{yz_rest},\eqref{y_dom},\eqref{z_dom},\nonumber\\
&&&\sum_{\pi\in \mathcal{P}_c}v_{\pi c}\leq 1,&& \hspace{-1cm}c\in C,\label{one_rest}\\
&&&\sum_{\pi\in \mathcal{P}_c} h_{e}^{\pi c}v_{\pi c}\leq y_{e},&&\hspace{-1cm} e\in E, \,c\in C,\label{hy_edge}\\
&&&\mbox{subtour elimination/connection constraints }&&\nonumber\\
&&&\mbox{ for the resulting structure},&&\nonumber\\
&&&v_{\pi c}\in\{0,1\},&&\hspace{-1cm}\pi\in \mathcal{P}_c,c\in C. \label{v_rest}
\end{alignat}

In the objective function, we maximize the total revenue in terms of travel time reduction of the resulting demand. As in $F1$, constraints \eqref{p_rest}--\eqref{yz_rest} ensure that only $p$ hubs are open, $p-1$ hub edges are operating, and that end-nodes of hub-edges are hubs, respectively. Constraints \eqref{one_rest} restrict that each commodity is transported either by, respectively, using the hub line or taking a direct path. In addition, constraints \eqref{hy_edge} ensure that commodities can only use paths whose hubs edges are opened. Constraints \eqref{v_rest} restrict $\bm{v}$-variables to be binary. 

Observe that, to obtain a valid formulation for \ED, we need to include a family of constraints establishing that the hub network is a path. Therefore, the formulation should avoid the existence of subtours or several connected components. In particular,  in ($F1_L$), $\bm f$-variables defined in \eqref{flowvar} and flow constraints (\eqref{flow1}-\eqref{flow3} and \eqref{f_dom}) could be used for that purpose.

Another way to ensure the path structure is to include the classic {\sl subtour elimination constraints} (SEC) \citep[][]{DeSa2015}. These constraints are usually applied and they could be expressed as follows:

\begin{equation}
\sum_{\stackrel{[i,m]\in E}{i,m\in S}}y_{[i,m]}\leq\sum_{i\in S\setminus\{s\}}z_i,\,\,\,S\subseteq N,s\in S.  \label{sec} \tag{\text{SEC}}
\end{equation}

Observe that the number of SEC constraints is exponential with respect to the number of nodes. Therefore, we do not consider all these constraints. Instead, we start with the SEC relaxation model and then, we sequentially include those that do not hold in the solution process by using a branch-and-cut algorithm.

 Another formulation for the \ED with a linear objective function can be derived by considering ${\bm y'}$-variables defined by \eqref{y'var}
and using constraints \eqref{MTZ}. As a result, we obtain the formulation below:

\begin{alignat}{3}
(F2_L)\quad& \max&  \quad & \sum_{c\in C}\sum_{\pi\in \mathcal{P}_c} g'_{jc}v_{\pi c}\nonumber\\
&\mbox{s.t.}&&\eqref{p_rest},\eqref{z_dom},\eqref{yp_rest2}-\eqref{y_dom2},\eqref{one_rest},\eqref{v_rest},\nonumber\\
&&&\sum_{\pi \in \mathcal{P}_c}h_{[k,m]}^{\pi c}v_{\pi c}\leq y'_{km}+y'_{mk},&& [k,m]\in E,c\in C.\label{hy_edge_n}
\end{alignat}

Observe that we have included constraints \eqref{yp_rest2}--\eqref{y_dom2} that are part of $F2$. For formulation $F2_L$, MTZ constraints avoid subtours and ensure connectivity. Similarly to constraints \eqref{hy_edge}, constraints \eqref{hy_edge_n} ensure that commodities use only paths whose hubs arcs are opened.

Analogously, 
\begin{alignat}{3}
(F2_L')\quad& \max&  \quad & \sum_{c\in C}\sum_{\pi\in \mathcal{P}_c} g'_{jc}v_{\pi c}\nonumber\\
&\mbox{s.t.}&&\eqref{p_rest},\eqref{z_dom},\eqref{yp_rest2},\eqref{MTZ}-\eqref{yz_new},\eqref{one_rest},\eqref{v_rest},\eqref{hy_edge_n}.\nonumber
\end{alignat}

Note that, as in the MINLP formulations, valid inequalities \eqref{new} can be included in $(F2_L')$ to avoid symmetries. In the next section, we present some valid inequalities for the above MILP formulations.

\subsection{Valid inequalities}

Given a commodity $c\in C$, we distinguish in set $\mathcal{P}_c$ four types of paths using the hub line. They are denoted as follows.
\begin{itemize}
	\item {\bf (ODH$_c$)-paths}. These are the paths in which all the nodes are hubs. This is, $o_c$ and $d_c$ are hubs connected by other open hubs.
	\item {\bf (DH$_c$)-paths}. This set of paths is composed by those in which $o_c$ is not a hub, but $d_c$ is one.
	\item {\bf (OH$_c$)-paths}. This set of paths is composed by those in which $d_c$ is not a hub, but $o_c$ is one.
	\item {\bf (ODNH$_c$)-paths}. These are the paths in which neither $o_c$ nor $d_c$ are hubs.
\end{itemize}

\begin{figure}[htpb]
	\centering
	\includegraphics[width=0.8\linewidth]{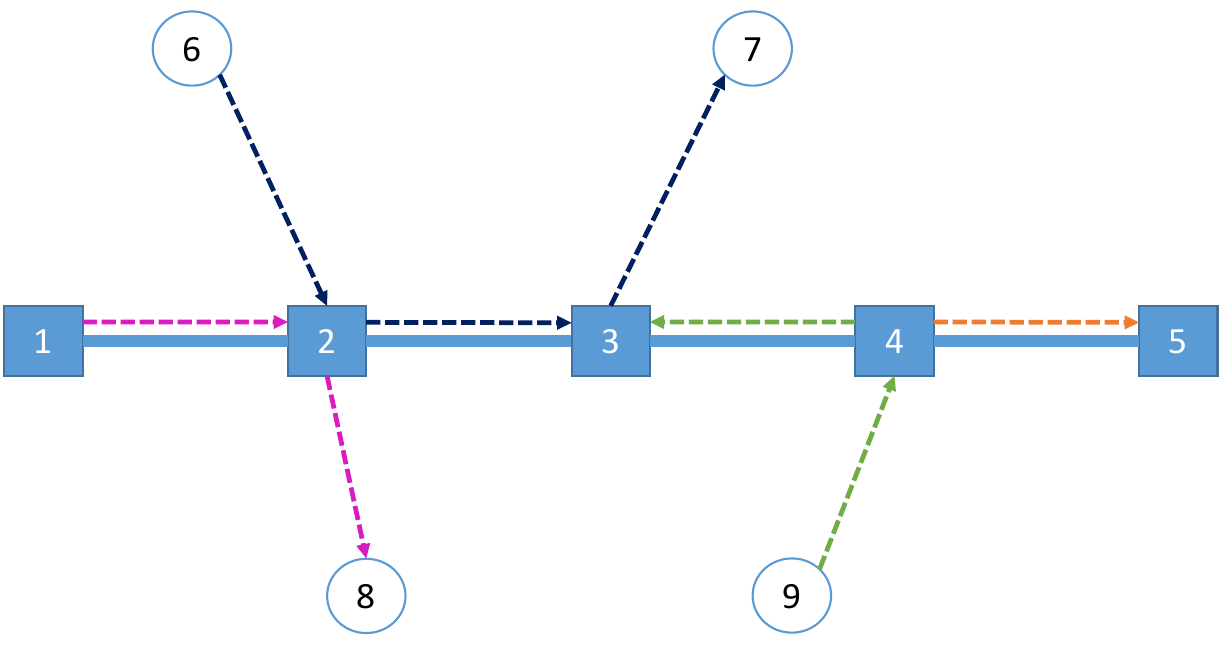}
	\caption{Example of the different types of path}
	\label{ejem}
\end{figure}

In Figure \ref{ejem} we can see an example of the different types of paths. Observe that in this case there is an already established hub line that traverses nodes 1, 2, 3, 4 and 5. The path between nodes 4 and 5 is a (ODH$_{(4,5)}$)-path since origin and destination are hubs. The path between 9 and 3 is a (DH$_{(9,3)}$)-path because the origin in not a hub but the destination is one. We can see an (OH$_{(1,8)}$)-path between nodes 1 and 8 since node 1 is a hub, but node 8 is not. Finally, between nodes 6 and 7 there is an (ODNH$_{(6,7)}$)-path since neither 6 nor 7 are hubs.

Attending to the different types of path that we have distinguished, we can include some valid inequalities for formulations $(F1_L)$, $(F2_L)$ and $(F2_L')$. The first family of valid inequalities is based on the following. Assume that, for a certain commodity $c$, all the candidates paths of type (ODH$_c$) and (DH$_c$) are known. Observe that in these types of paths, $d_c$ must be a hub. Consequently, we obtain the next constraints:

\begin{eqnarray}
\sum_{\pi\in \mbox{{\scriptsize (ODH$_c$)-paths}}}
v_{\pi c}+\sum_{\pi\in \mbox{{\scriptsize (DH$_c$)-paths}}}v_{\pi c}\leq z_{d_c},\quad c\in C. \label{desthub}
\end{eqnarray}

Similarly (ODH$_c$)-paths and (OH$_c$)-paths, for a commodity $c$, satisfy that $o_c$ is a hub. Then, 

\begin{eqnarray}
\sum_{\pi\in \mbox{{\scriptsize (ODH$_c$)-paths}}}v_{\pi c}+\sum_{\pi\in \mbox{{\scriptsize (OH$_c$)-paths}}}v_{\pi c}\leq z_{o_c},\quad c\in C. \label{orhub}
\end{eqnarray}
Recall that for formulations $(F1_L)$, $(F2_L)$ and $(F2'_L)$ to be valid, we need to know a priori all the candidate OD paths. Section \ref{paths} focuses on an efficient procedure to obtain the set of such OD paths.

\subsection{Generating all candidate OD paths}\label{paths}

The MILP formulations ($F1_L$, $F2_L$ and $F2_L'$) proposed in the previous section use the candidate OD paths of each commodity as input. Consequently, it is necessary to develop a procedure to obtain all the candidate paths. With this aim, we implement an algorithm based on the use of function \emph{all\_simple\_paths} that can be found in the \emph{Networkx} library of  Python. This function allows us to obtain all simple paths between a pair of nodes ($o_c$ and $d_c$) with a certain length. 

The general scheme to find all candidate paths of each $c\in C$ is the following.

\noindent\underline{{\bf Step 0.}} Let consider a graph  $G=(N,E)$, where $N$ is the set of demand points and potential hubs, and $E$ is the set of existing edges between them. Moreover, $t_{ij}$ represents the time needed to traverse edge $e=[i,j]\in E$.

\noindent\underline{{\bf Step 1.}}  Given a commodity $c\in C$, build the {auxiliary} graph $G_c=(N_c,A_c)$ and travel times $t^c$ by using Algorithm \ref{alg-1}.

\noindent\underline{{\bf Step 2.}} Given  $G_c=(N_c,A_c)$ and $t^c$, search for all the possible paths with a number of arcs smaller than or equal to $p+1$ (the number of open hubs is $p$). To carry out this step, function  {\sl all\_simple\_paths} is used.

\noindent\underline{{\bf Step 3 (discarding paths).}} Discard a path whenever the time to traverse it is greater than or equal to $t_{o_c d_c}$. Thus, the final paths are only those with a time smaller than $t_{o_cd_c}$. 

Given that an arc $(i,j)\in A_c$ is part of a certain path (denoted as $\pi$), if the path $[o_c, i, j, d_c]$ has a travel time ($t_{o_ci}+\tilde{t}^a_i+\alpha t_{ij}+\tilde{t}^e_j+t_{jd_c}$) smaller than or equal to the time of $\pi$, then $\pi$ is not included in the final set of paths of commodity $c$. This is due to the fact that, if we open the hub edge $[i,j]$, commodity $c$ will prefer to route its demand only using those hubs than using $\pi$ since the associated travel time is smaller. 

After applying this procedure, the set of candidate paths for commodity $c$, $\mathcal{P}_c$, is determined. Algorithm \ref{alg1} details the procedure to obtain all possible paths in a hub line composed by $p$ hubs.

\begin{algorithm}[H]
	\label{alg1}
	\KwData{Graph $G=(N,E)$ with $t_{ij}$ for $e=[i,j]\in E$ and $c\in C$.} 
	\KwResult{All possible (ODH$_c$, DH$_c$, OH$_c$ and ODNH$_c$)-paths for $c\in C$.} 
	Build the {auxiliary} graph $G_c=(N_c,A_c)$ and travel times $t^c$ by using Algorithm \ref{alg-1}.\\
		Search for the paths, $\pi_c$, associated to $c$ (all\_simple\_paths($G_c,o_c,d_c,p+1$)) such that:\\
		$$\#_{\mbox{arcs}}(\pi_c)\leq p+1 \mbox{ and time}(\pi_c)\leq t_{o_c d_c}$$
		.\\
		Where $\#_{\mbox{arcs}}(\pi_c)$ is the number of arcs in $\pi_c$ and time$(\pi_c)$ is the time of $\pi_c$. \\
		
		 Once a path $\pi_c$ is obtained,\\
		\If{there exists an arc $(i,j)$ of $\pi_c$ such that $t_{o_ci}+\tilde{t}^a_i+\alpha t_{ij}+\tilde{t}^e_j+t_{jd_c}\leq $time$(\pi_c)$ }{
		remove $\pi_c$ from the set of paths associated with $(o_c,d_c)$.}

	Given a candidate path $\pi_c$, identify type of $\pi_c$:\\
		\If{$i'= o_c$ and  $j'= d_c$ }{ 
		$\pi_c$ is type ODH$_c$
		}
		\If{$i'= o_c$ and  $j'\neq d_c$ }{ 
		$\pi_c$ is type OH$_c$
		}
		\If{$i'\neq o_c$ and  $j'= d_c$ }{ 
		$\pi_c$ is type DH$_c$
		}	
		\If{$i'\neq o_c$ and  $j'\neq d_c$ }{ 
		$\pi_c$ is type ODNH$_c$
		}

	\Return All possible paths for $c$.
	\caption{Obtaining all (ODH$_c$, DH$_c$, OH$_c$ and ODNH$_c$)-paths for a commodity $c\in C$.}
\end{algorithm}

\section{Computational Results}
\label{computational}

 We next carry out two main experimental studies. In Section \ref{cab}, we first compare the computational results provided by the previously described formulations for \ED. For these experiments, the benchmark dataset of Civil Aeronautics Board (CAB) was used \citep[see,][]{Okelly}. The main goal of this study is to determine which of the formulations (MINLP or MILP) provides the best performance. In Section \ref{subsecMontreal}, we then analyze the quality of the solutions provided by \ED in a case study we generated for the city of Montreal. We used real data provided by The Autorité Régionale de Transport Métropolitain (ARTM). We compare the performance of \ED using different values of the parameters involved in the \ED, and we compare the results with the solutions obtained by the \HLLP.

\subsection{Comparison of formulations}\label{cab}

This section is devoted to the comparison of the computational results provided by the described formulations for \ED.  For this set of experiments, we used the CAB dataset. 
We note that these instances have a symmetric OD demand matrices. We considered instances with $n=10, 15, 20, 25$. Our computational studies were performed on an Intel(R) Xeon(R) W-2245 CPU 256 GB RAM computer.

For the nonlinear models $(F1)$, $(F2)$ and $(F2')$, the global optimization software BARON 21.1.13 was used through the AMPL modeling language. To study the performance of $(F1)$, $(F2)$ and $(F2')$, CAB instances with complete graphs of 10 and 15 nodes were examined. In addition, regarding the parameters of the problem, we use $R_c=(1+\gamma_c)t_{o_c d_c}$ selecting  $\gamma_c \in[0,1]$ randomly, $r=1.7$, $\alpha\in\{0.2,0.5,0.8\}$ and $p\in\{2,3,5\}$ \citep[see,][]{FKelly1989,Zetina2019}. Similarly to \cite{DeSa2015}, we consider that the access and exit times do not depend on the node, i.e., $\tilde{t}^a_k=\tilde{t}^a$ and $\tilde{t}^e_k=\tilde{t}^e$ for all $k\in N$. In particular, we define access and exit time as a proportion of the average travel time:
$$\tilde{t}^a=\tilde{t}^e=\displaystyle \vartheta\frac{\sum_{(i,j)\in A}t_{ij}}{n\cdot (n-1)},$$ 
where $\vartheta$ is fixed to $0.1$ for these computational results. Table \ref{tab_nonlinear} reports the results of the developed nonlinear formulations for \ED using a time limit of two hours. The first column presents the number of nodes of the considered instances ($n$), the second column reports the number of open hubs ($p$) and, then the used values of the parameter $\alpha$ are specified. Recall that BARON requires to include adequate upper and lower bounds for the nonlinear terms appearing in the objective function of the models. Thus, we used the preprocessing phase explained in Section \ref{impro_nolinear}, whose running time in seconds is shown in column  $t_{prep}$. 

The remaining columns report results for formulations $(F1)+\eqref{lb_ub}-\eqref{fixn}$, $(F2)+\eqref{lb_ub}-\eqref{fixn}$, $(F2')+\eqref{lb_ub}-\eqref{fixn}$ and $(F2')+\eqref{new}-\eqref{fixn}$. Observe that, in all formulations, we have added fixing criteria and valid inequalities proposed in Subsection \ref{impro_nolinear}. We have included them since preliminary results show that \eqref{fix0}--\eqref{fixn} improve the computational performance of the formulations. For each formulation, we show two columns. The first one is the gap if the instance is not solved in two hours, and the second column is the solving time. Note that column ``Gap'' is calculated as: $\frac{UB-LB}{LB}\times 100$, where $UB$ is the best-obtained upper bound in 7,200 seconds and $LB$ is the best-obtained solution in this time limit. 

Regarding the results of Table \ref{tab_nonlinear}, we can assure that the solving times become bigger when $n$ and $p$ parameter values increase and the $\alpha$ value decreases. We observe that all formulations provide a similar time performance.  For the biggest considered instances ($n=15, p=5$), the best results are provided by $(F2)$ which is able to solve the case with $\alpha=0.8$ in less than $500$ seconds and provides the best gap for $\alpha=0.5$ and $\alpha=0.2$.

\begin{table}[htbp]
  \centering
  \caption{Time results of formulations $(F1)+\eqref{lb_ub}-\eqref{fixn}$, $(F2)+\eqref{lb_ub}-\eqref{fixn}$, $(F2')+\eqref{lb_ub}-\eqref{fixn}$ and $(F2')+\eqref{new}-\eqref{fixn}$ for $n=10$ and $n=15$}
  	\begin{adjustbox}{max width=1\textwidth}
    \begin{tabular}{ccrrrrrrrrrr}
    \toprule
    \multirow{2}[2]{*}{$n$} & \multirow{2}[2]{*}{$p$} & \multicolumn{1}{c}{\multirow{2}[2]{*}{$\alpha$}} & \multicolumn{1}{c}{\multirow{2}[2]{*}{$t_{prep}$}} & \multicolumn{2}{c}{ $(F1)+\eqref{lb_ub}-\eqref{fixn}$} & \multicolumn{2}{c}{$(F2)+\eqref{lb_ub}-\eqref{fixn}$} & \multicolumn{2}{c}{$(F2')+\eqref{lb_ub}-\eqref{fixn}$} & \multicolumn{2}{c}{ $(F2')+\eqref{new}-\eqref{fixn}$} \\
    \cmidrule(lr){5-6}
   \cmidrule(lr){7-8}
   \cmidrule(lr){9-10}
    \cmidrule(lr){11-12}
          &       &       &       & \multicolumn{1}{c}{Gap} & \multicolumn{1}{c}{Time} & \multicolumn{1}{c}{Gap} & \multicolumn{1}{c}{Time} & \multicolumn{1}{c}{Gap} & \multicolumn{1}{c}{Time} & \multicolumn{1}{c}{Gap} & \multicolumn{1}{c}{Time} \\
    \midrule
    
    \multirow{3}[1]{*}{10} & \multirow{3}[1]{*}{2} & 0.8   & 1.97  & 0.00  & 27.14 & 0.00  & 26.13 & 0.00  & \textbf{25.58} & 0.00  & 27.06 \\
          &       & 0.5   & 2.08  & 0.00  & 82.91 & 0.00  & \textbf{73.81} & 0.00  & 105.59 & 0.00  & 82.39 \\
          &       & 0.2   & 2.14  & 0.00  & 343.59 & 0.00  & \textbf{273.31} & 0.00  & 274.19 & 0.00  & 386.78 \\
    \midrule
    \multirow{3}[1]{*}{10} & \multirow{3}[1]{*}{3} & 0.8   & 1.96  & 0.00  & 76.45 & 0.00  & \textbf{40.72} & 0.00  & 57.48 & 0.00  & 43.36 \\
          &       & 0.5   & 2.07  & 0.00  & 130.63 & 0.00  & \textbf{70.56} & 0.00  & 304.45 & 0.00  & 112.59 \\
          &       & 0.2   & 2.15  & 0.00  & 1304.59 & 0.00  & \textbf{912.56} & 0.00  & 2891.64 & 0.00  & 3619.63 \\
          \midrule
    \multirow{3}[1]{*}{10} & \multirow{3}[1]{*}{5} & 0.8   & 1.98  & 0.00  & \textbf{37.78} & 0.00  & 72.09 & 0.00  & 51.34 & 0.00  & 52.11 \\
          &       & 0.5   & 2.07  & 0.00  & 525.44 & 0.00  & 1065.69 & 0.00  & 334.98 & 0.00  & \textbf{91.33} \\
          &       & 0.2   & 2.16  & 17.97 & 7200.00 & 21.61 & 7200.00 & 24.21 & 7200.00 & 19.21 & 7200.00 \\
    \midrule
    \multirow{3}[1]{*}{15} & \multirow{3}[1]{*}{2} & 0.8   & 2.91  & 0.00  & 1184.73 & 64.22 & 7200.00 & 0.00  & 1314.39 & 0.00  & \textbf{1176.84} \\
          &       & 0.5   & 3.30  & 0.00  & 2287.27 & 0.00  & 2929.16 & 0.00  & 2190.66 & 0.00  & \textbf{475.78} \\
          &       & 0.2   & 3.49  & 0.00  & 2248.66 & 38.45 & 7200.00 & 41.65 & 7200.00 & 0.00  & \textbf{2219.20} \\
          \midrule
    \multirow{3}[1]{*}{15} & \multirow{3}[1]{*}{3} & 0.8   & 2.96  & 0.00  & 1373.92 & 0.00  & \textbf{609.13} & 0.00  & 2176.61 & 0.00  & 2042.70 \\
          &       & 0.5   & 3.29  & 0.00  & \textbf{4279.56} & 16.45 & 7200.00 & 10.62 & 7200.00 & 15.68 & 7200.00 \\
          &       & 0.2   & 3.48  & 34.68 & 7200.00 & 34.68 & 7200.00 & 95.67 & 7200.00 & 34.68 & 7200.00 \\
    \midrule
    \multirow{3}[2]{*}{15} & \multirow{3}[2]{*}{5} & 0.8   & 3.00  & 0.00  & 1534.98 & 0.00  & \textbf{425.97} & 0.00  & 1164.50 & 0.00  & 3371.55 \\
          &       & 0.5   & 3.28  & 23.50 & 7200.00 & 16.27 & 7200.00 & 16.27 & 7200.00 & 22.62 & 7200.00 \\
          &       & 0.2   & 3.51  & 37.82 & 7200.00 & 36.89 & 7200.00 & 88.29 & 7200.00 & 39.98 & 7200.00 \\
    \bottomrule
    \end{tabular}%
    \end{adjustbox}
  \label{tab_nonlinear}%
\end{table}%

\begin{table}[htbp]
  \centering
  \caption{Time results of formulations $(F1_L)$+\eqref{flow1}-\eqref{flow3}, $(F1_L)$+\eqref{sec}, $(F2_L)$, $(F2'_L)$ and $(F2'_L)$+\eqref{new} for $n=10$ and $n=15$. In all cases, we include valid inequalities \eqref{desthub} and \eqref{orhub}.}
  	\begin{adjustbox}{max width=1\textwidth}
    \begin{tabular}{cclrrrrrrrrrrrr}
    \toprule
    \multirow{2}[2]{*}{$n$} & \multirow{2}[2]{*}{$p$} & \multicolumn{1}{c}{\multirow{2}[2]{*}{$\alpha$}} & \multicolumn{1}{c}{\multirow{2}[2]{*}{$n_{path}$}} & \multicolumn{1}{c}{\multirow{2}[2]{*}{$t_{path}$}} & \multicolumn{2}{c}{$(F1_L)$+\eqref{flow1}-\eqref{flow3}} & \multicolumn{2}{c}{$(F1_L)$+\eqref{sec}} & \multicolumn{2}{c}{$(F2_L)$} & \multicolumn{2}{c}{$(F2'_L)$} & \multicolumn{2}{c}{$(F2'_L)$+\eqref{new}} \\
     \cmidrule(lr){6-7}
   \cmidrule(lr){8-9}
   \cmidrule(lr){10-11}
    \cmidrule(lr){12-13}\cmidrule(lr){14-15}
          &       &       &       &       & \multicolumn{1}{c}{Gap$_{lp}$} & \multicolumn{1}{c}{$t$} & \multicolumn{1}{c}{Gap$_{lp}$} & \multicolumn{1}{c}{$t$} & \multicolumn{1}{c}{Gap$_{lp}$} & \multicolumn{1}{c}{$t$} & \multicolumn{1}{c}{Gap$_{lp}$} & \multicolumn{1}{c}{$t$} & \multicolumn{1}{c}{Gap$_{lp}$} & \multicolumn{1}{c}{$t$} \\
    \midrule
    \multirow{3}[2]{*}{10} & \multirow{3}[2]{*}{2} & 0.8   & 80    & 0.64  & 0.00  & 0.09  & 0.00  & \textbf{0.06} & 0.00  & 0.08  & 0.00  & 0.06  & 0.00  & 0.08 \\
          &       & 0.5   & 338   & 0.71  & 0.00  & 0.14  & 0.00  & \textbf{0.11} & 0.00  & \textbf{0.11} & 0.00  & 0.14  & 0.00  & 0.14 \\
          &       & 0.2   & 674   & 1.01  & 0.00  & 0.20  & 0.00  & \textbf{0.13} & 0.00  & \textbf{0.13} & 0.00  & 0.13  & 0.00  & 0.16 \\
    \midrule
    \multirow{3}[2]{*}{10} & \multirow{3}[2]{*}{3} & 0.8   & 190   & 0.66  & 0.00  & 0.08  & 0.00  & \textbf{0.06} & 0.00  & \textbf{0.06} & 0.00  & 0.08  & 0.00  & \textbf{0.06} \\
          &       & 0.5   & 1134  & 0.77  & 0.00  & 0.16  & 0.00  & 0.17  & 0.00  & 0.14  & 0.00  & \textbf{0.13} & 0.00  & 0.14 \\
          &       & 0.2   & 2980  & 0.90  & 0.00  & 0.24  & 0.00  & 0.22  & 0.00  & \textbf{0.20} & 0.00  & 0.20  & 0.00  & 0.31 \\
    \midrule
    \multirow{3}[2]{*}{10} & \multirow{3}[2]{*}{5} & 0.8   & 292   & 1.02  & 6.55  & 0.22  & 6.55  & \textbf{0.11} & 6.55  & 0.19  & 6.55  & 0.17  & 6.55  & 0.19 \\
          &       & 0.5   & 3836  & 1.55  & 5.53  & 0.44  & 5.53  & 2.23  & 5.53  & \textbf{0.41} & 5.53  & 0.48  & 5.53  & 0.53 \\
          &       & 0.2   & 32554 & 4.94  & 4.35  & 4.16  & 4.35  & 6.84  & 4.35  & \textbf{4.12} & 4.35  & 8.03  & 4.35  & 11.79 \\
    \midrule
    \multirow{3}[2]{*}{15} & \multirow{3}[2]{*}{2} & 0.8   & 214   & 0.92  & 0.00  & 0.31  & 0.00  & 0.23  & 0.00  & \textbf{0.23} & 0.00  & 0.25  & 0.00  & 0.25 \\
          &       & 0.5   & 1202  & 1.23  & 0.00  & 0.61  & 0.00  & \textbf{0.47} & 0.00  & 0.58  & 0.00  & 0.67  & 0.00  & 0.83 \\
          &       & 0.2   & 2372  & 1.33  & 0.00  & 0.87  & 0.00  & \textbf{0.56} & 0.00  & 0.69  & 0.00  & 0.67  & 0.00  & 1.05 \\
    \midrule
    \multirow{3}[2]{*}{15} & \multirow{3}[2]{*}{3} & 0.8   & 756   & 1.05  & 0.00  & 0.28  & 0.00  & \textbf{0.25} & 0.00  & 0.33  & 0.00  & 0.28  & 0.00  & 0.31 \\
          &       & 0.5   & 4822  & 1.64  & 0.00  & 0.75  & 0.00  & \textbf{0.67} & 0.00  & 0.93  & 0.00  & 0.72  & 0.00  & 0.75 \\
          &       & 0.2   & 14438 & 2.46  & 0.00  & 1.34  & 0.00  & \textbf{1.22} & 0.00  & 1.45  & 0.00  & 1.42  & 0.00  & 1.44 \\
    \midrule
    \multirow{3}[2]{*}{15} & \multirow{3}[2]{*}{5} & 0.8   & 2028  & 11.06 & 0.00  & 0.34  & 0.00  & \textbf{0.31} & 0.00  & 0.42  & 0.00  & 0.36  & 0.00  & 0.37 \\
          &       & 0.5   & 31010 & 21.55 & 1.84  & \textbf{3.54} & 1.84  & 9.30  & 1.84  & 5.59  & 1.84  & 26.96 & 1.84  & 14.14 \\
          &       & 0.2   & 414430 & 73.76 & 10.08 & 859.46 & 10.08 & 408.98 & 10.08 & 399.07 & 10.08 & 219.34 & 10.08 & \textbf{202.20} \\
    \bottomrule
    \end{tabular}%
    \end{adjustbox}
  \label{tab10_15}%
\end{table}%

In general, these MINLP formulations require big solving times for instances of $n=10$ and $n=15$. This limited performance makes necessary the use of alternative formulations to solve \ED. Thus, we focus on the results provided by  the MILP formulations proposed in Section \ref{linear}. They were implemented in Python and solved using CPLEX 20.1.0.

Table \ref{tab10_15} shows the results for formulations $(F1_L)$+\eqref{flow1}-\eqref{flow3}, $(F1_L)$+\eqref{sec}, $(F2_L)$, $(F2_L')$ and $(F2_L')$+\eqref{new}
using CAB instances with $n=10$ and $n=15$. First column reports the number of nodes of the instances ($n$), the second column specifies the number of open hubs ($p$), and the third column reports the $\alpha$ values. Observe that, for using proposed MILP formulations for \ED, it is necessary to know the set of possible paths associated with each commodity. In order to obtain these sets, the algorithm described in Subsection \ref{paths} was implemented.

Observe that the described procedure in Algorithm \ref{alg1} can be carried out independently for each commodity $c$. Regarding the implementation of this path search, we can think about a main for-loop which goes through all commodities and applies the search for each of them. However, some execution time can be reduced if we try to run the procedure for each commodity in parallel. 

In order to apply the procedure for different commodities in a simultaneous way, we use {\sl concurrent.futures} model of Python, which provides a high-level interface for asynchronously executing callable. Particularly, we apply {\sl ProcessPoolExecutor} class that allows to display different subprocesses at the same time. The maximum number of simultaneous sub-processes is fixed to the number of CPUs divided by two. More details about this class can be found in \url{https://docs.python.org/3/library/concurrent.futures.html\#concurrent.futures.ProcessPoolExecutor}

The number of obtained paths using Algorithm \ref{alg1} is reported in column $n_{path}$. Besides, column $t_{path}$ reports the necessary running time to obtain the path candidates. The remaining columns report the gap with the linear relaxation and the total running time of the formulations ($(F1_L)$+\eqref{flow1}-\eqref{flow3}, $(F1_L)$+\eqref{sec}, $(F2_L)$, $(F2'_L)$ and $(F2'_L)+\eqref{new}$). In these formulations we have included valid inequalities \eqref{desthub} and \eqref{orhub} since they show to provide better time results in a preliminary study. In this table, in contrast to Table \ref{tab_nonlinear}, the gaps at termination are not reported since all instances were solved within the time limit. 

In general, if we compare the results of these formulations with respect nonlinear formulations, we can claim that MILP formulations result to be more efficient.  Observe that the running times of formulations in Table \ref{tab10_15} are smaller, and all the instances can be solved in less than 203 seconds. These results show the advantages of using formulations developed in Section \ref{linear}. Particularly, we can see in boldface the formulation that provides the best time result for each instance. We can remark that the best formulation for these instance sizes seems to be $(F1_L)$+\eqref{sec}. However, for $n=15$, $p=5$ and $\alpha=0.2$ the best time result is provided by $(F2_L')$+\eqref{new}.

\begin{sidewaystable}[htbp]
  \centering
  \caption{Time results of formulations $(F1_L)$+\eqref{flow1}-\eqref{flow3}, $(F1_L)$+\eqref{sec}, $(F2_L)$,$(F2'_L)$ and $(F2'_L)+\eqref{new}$,  with valid inequalities \eqref{desthub}+\eqref{orhub} for $n=20$}
  \begin{adjustbox}{max width=1\textwidth}
    \begin{tabular}{cclrrrrrrrrrrrrrrrrr}
    \toprule
    \multirow{2}[2]{*}{$p$} & \multirow{2}[2]{*}{$\alpha$} & \multicolumn{1}{c}{\multirow{2}[2]{*}{(\% edges)}} & \multicolumn{1}{c}{\multirow{2}[2]{*}{$n_{path}$}} & \multicolumn{1}{c}{\multirow{2}[2]{*}{$t_{path}$}} & \multicolumn{3}{c}{$(F1_L)$+\eqref{flow1}-\eqref{flow3}} & \multicolumn{3}{c}{$(F1_L)$+\eqref{sec}} & \multicolumn{3}{c}{$(F2_L)$} & \multicolumn{3}{c}{$(F2'_L)$} & \multicolumn{3}{c}{$(F2'_L)$+\eqref{new}} \\
      \cmidrule(lr){6-8}
      \cmidrule(lr){9-11}
      \cmidrule(lr){12-14}
      \cmidrule(lr){15-17}
      \cmidrule(lr){18-20}
          &       &       &       &       & \multicolumn{1}{c}{Gap$_{lp}$} & \multicolumn{1}{c}{Gap} & \multicolumn{1}{c}{$t$} & \multicolumn{1}{c}{Gap$_{lp}$} & \multicolumn{1}{c}{Gap} & \multicolumn{1}{c}{$t$} & \multicolumn{1}{c}{Gap$_{lp}$} & \multicolumn{1}{c}{Gap} & \multicolumn{1}{c}{$t$} & \multicolumn{1}{c}{Gap$_{lp}$} & \multicolumn{1}{c}{Gap} & \multicolumn{1}{c}{$t$} & \multicolumn{1}{c}{Gap$_{lp}$} & \multicolumn{1}{c}{Gap} & \multicolumn{1}{c}{$t$} \\
       \midrule
    \multirow{3}[2]{*}{3} & \multirow{3}[2]{*}{0.8} & 0.4   & 766   & 1.49  & 0.00  & 0.00  & 0.66  & 0.00  & 0.00  & \textbf{0.62} & 0.00  & 0.00  & 0.74  & 0.00  & 0.00  & 0.74  & 0.00  & 0.00  & 18.88 \\
          &       & 0.5   & 1148  & 1.61  & 0.00  & 0.00  & 0.88  & 0.00  & 0.00  & \textbf{0.68} & 0.00  & 0.00  & 0.89  & 0.00  & 0.00  & 0.85  & 0.00  & 0.00  & 0.86 \\
          &       & 0.6   & 1550  & 1.76  & 0.00  & 0.00  & 0.82  & 0.00  & 0.00  & \textbf{0.76} & 0.00  & 0.00  & 0.92  & 0.00  & 0.00  & 0.90  & 0.00  & 0.00  & 0.92 \\
    \midrule
    \multirow{3}[2]{*}{3} & \multirow{3}[2]{*}{0.5} & 0.4   & 4704  & 2.46  & 0.00  & 0.00  & 1.62  & 0.00  & 0.00  & \textbf{1.51} & 0.00  & 0.00  & 1.77  & 0.00  & 0.00  & 1.77  & 0.00  & 0.00  & 1.83 \\
          &       & 0.5   & 7070  & 2.76  & 0.00  & 0.00  & 1.97  & 0.00  & 0.00  & \textbf{1.77} & 0.00  & 0.00  & 2.16  & 0.00  & 0.00  & 2.09  & 0.00  & 0.00  & 2.13 \\
          &       & 0.6   & 9543  & 3.09  & 0.00  & 0.00  & 2.25  & 0.00  & 0.00  & \textbf{2.03} & 0.00  & 0.00  & 2.42  & 0.00  & 0.00  & 2.40  & 0.00  & 0.00  & 2.32 \\
    \midrule
    \multirow{3}[2]{*}{3} & \multirow{3}[2]{*}{0.2} & 0.4   & 13054 & 3.34  & 0.00  & 0.00  & 2.58  & 0.00  & 0.00  & \textbf{2.42} & 0.00  & 0.00  & 2.85  & 0.00  & 0.00  & 2.94  & 0.00  & 0.00  & 2.97 \\
          &       & 0.5   & 19722 & 3.92  & 0.00  & 0.00  & 3.28  & 0.00  & 0.00  & \textbf{3.05} & 0.00  & 0.00  & 3.64  & 0.00  & 0.00  & 3.56  & 0.00  & 0.00  & 3.60 \\
          &       & 0.6   & 27379 & 4.64  & 0.00  & 0.00  & 4.01  & 0.00  & 0.00  & \textbf{3.79} & 0.00  & 0.00  & 4.33  & 0.00  & 0.00  & 4.33  & 0.00  & 0.00  & 4.31 \\
    \midrule
    \multirow{3}[2]{*}{5} & \multirow{3}[2]{*}{0.8} & 0.4   & 1195  & 4.31  & 0.99  & 0.00  & 0.83  & 1.02  & 0.00  & \textbf{0.78} & 1.02  & 0.00  & 1.09  & 1.02  & 0.00  & 1.08  & 0.93  & 0.00  & 1.10 \\
          &       & 0.5   & 2142  & 8.46  & 2.36  & 0.00  & \textbf{1.19} & 2.36  & 0.00  & 1.86  & 2.36  & 0.00  & 1.91  & 2.36  & 0.00  & 1.60  & 2.36  & 0.00  & 1.61 \\
          &       & 0.6   & 3384  & 16.41 & 3.11  & 0.00  & \textbf{1.32} & 3.11  & 0.00  & 2.59  & 3.11  & 0.00  & 3.26  & 3.11  & 0.00  & 2.72  & 3.11  & 0.00  & 2.65 \\
    \midrule
    \multirow{3}[2]{*}{5} & \multirow{3}[2]{*}{0.5} & 0.4   & 10636 & 7.43  & 0.99  & 0.00  & \textbf{2.87} & 0.99  & 0.00  & 5.07  & 0.99  & 0.00  & 8.74  & 0.99  & 0.00  & 10.64 & 0.99  & 0.00  & 16.25 \\
          &       & 0.5   & 22966 & 14.94 & 1.10  & 0.00  & \textbf{4.21} & 1.10  & 0.00  & 7.68  & 1.10  & 0.00  & 12.34 & 1.10  & 0.00  & 13.57 & 1.10  & 0.00  & 7.40 \\
          &       & 0.6   & 42385 & 28.56 & 0.58  & 0.00  & 8.91  & 0.58  & 0.00  & 9.56  & 0.58  & 0.00  & 12.42 & 0.58  & 0.00  & 14.28 & 0.58  & 0.00  & \textbf{7.42} \\
    \midrule
    \multirow{3}[2]{*}{5} & \multirow{3}[2]{*}{0.2} & 0.4   & 115494 & 23.52 & 0.13  & 0.00  & 13.94 & 0.13  & 0.00  & \textbf{13.66} & 0.13  & 0.00  & 24.94 & 0.13  & 0.00  & 33.04 & 0.13  & 0.00  & 19.55 \\
          &       & 0.5   & 270086 & 53.18 & 1.16  & 0.00  & \textbf{77.23} & 1.16  & 0.00  & 110.67 & 1.16  & 0.00  & 97.88 & 1.16  & 0.00  & 265.18 & 1.16  & 0.00  & 131.09 \\
          &       & 0.6   & 556068 & 105.55 & 0.98  & 0.00  & 126.88 & 0.98  & 0.00  & \textbf{116.88} & 0.98  & 0.00  & 189.05 & 0.98  & 0.00  & 462.13 & 0.98  & 0.00  & 172.98 \\
    \midrule
    \multirow{3}[2]{*}{7} & \multirow{3}[2]{*}{0.8} & 0.4   & 1206  & 77.20 & 3.39  & 0.00  & \textbf{1.11} & 3.50  & 0.00  & 3.01  & 3.50  & 0.00  & 3.18  & 3.39  & 0.00  & 3.79  & 3.29  & 0.00  & 1.68 \\
          &       & 0.5   & 2184  & 291.38 & 4.31  & 0.00  & \textbf{1.26} & 4.31  & 0.00  & 4.72  & 4.31  & 0.00  & 6.87  & 4.31  & 0.00  & 2.76  & 4.31  & 0.00  & 3.52 \\
          &       & 0.6   & 3505  & 957.73 & 6.29  & 0.00  & \textbf{2.42} & 6.29  & 0.00  & 10.01 & 6.29  & 0.00  & 17.22 & 5.82  & 0.00  & 10.12 & 6.29  & 0.00  & 9.79 \\
    \midrule
    \multirow{3}[2]{*}{7} & \multirow{3}[2]{*}{0.5} & 0.4   & 14023 & 123.00 & 3.24  & 0.00  & \textbf{15.52} & 3.24  & 0.00  & 24.65 & 3.24  & 0.00  & 48.27 & 3.24  & 0.00  & 31.49 & 3.20  & 0.00  & 50.38 \\
          &       & 0.5   & 38976 & 458.00 & 3.73  & 0.00  & \textbf{29.18} & 3.73  & 0.00  & 44.68 & 3.73  & 0.00  & 51.50 & 3.27  & 0.00  & 41.83 & 3.71  & 0.00  & 57.93 \\
          &       & 0.6   & 91821 & 1441.78 & 4.08  & 0.00  & \textbf{52.81} & 4.08  & 0.00  & 67.61 & 4.08  & 0.00  & 145.44 & 3.69  & 0.00  & 116.36 & 4.08  & 0.00  & 132.94 \\
    \midrule
    \multirow{3}[2]{*}{7} & \multirow{3}[2]{*}{0.2} & 0.4   & 760917 & 369.67 & 6.71  & 0.00  & \textbf{663.19} & 6.71  & 0.00  & 677.15 & 6.71  & 0.00  & 1761.22 & 6.54  & 0.00  & 1258.53 & 6.71  & 0.01  & 1479.84 \\
          &       & 0.5   & 2894947 & 1370.86 & 8.41  & (3)3.97  & \textbf{5387.17} & 9.26  & (3)4.88  & 5777.46 & 10.25 & (4)8.17  & 6768.62 & 8.54  & (5)5.97  & 7519.48 & 8.53  & (5)6.94  & 7525.06 \\
          &       & 0.6   & 8931072 & 4144.62 & 8.75  & (5)7.08  & 8195.44 & 9.87  & (4)8.05  & \textbf{8032.42} & 8.96  & (5)8.00  & 8206.01 & 8.75  & (4)7.19  & 8099.77 & 9.97  & (5)9.00  & 8211.58 \\
    \bottomrule
    \end{tabular}%
    \end{adjustbox}
  \label{tab_20}%
\end{sidewaystable}%

In Table \ref{tab_20}, we report more results for the MILP formulations proposed in Section \ref{linear}. In this case, we consider $n=20$, $p\in\{3,5,7\}$ and $\alpha\in\{0.2,0.5,0.8\}$. We did not use the complete graph provided by CAB instances for these parameters. Instead, we created more sparse graphs by only considering a subset of edges as possible hub edges. Particularly, we sorted the edges in non-decreasing order with respect to their associated times, and we discarded the edges with a time among the $10\%$ smallest ones or among the $10\%$ largest ones. Then, we selected randomly $40\%$, $50\%$ or $60\%$  among the remaining edges. The parameter that determines this  sparsity is denoted as (\% edges), and it is detailed in the third column of Table \ref{tab_20}. For each combination of these parameter values, the corresponding row of Table \ref{tab_20} shows the average results for five randomly generated instances.

Similarly to Table \ref{tab10_15}, Table \ref{tab_20} reports in columns $n_{path}$ and $t_{path}$ the average number of obtained paths and the average running times of Algorithm \ref{alg1}.
It should be remarked that this algorithm allows to generate a big number of paths in a reasonable time.
For instance, for $n=20$, $p=7$, $\alpha=0.2$, (\%edges)$=0.6$, the average time is 4,144.62, and the average number of generated paths is 8,931,072.

The remaining columns report the LP gap, the gap between the best solution and the best upper bound within the time limit, and the running times of formulations $(F1_L)$+\eqref{flow1}-\eqref{flow3}, $(F1_L)$+\eqref{sec}, $(F2_L)$, $(F2_L')$  and $(F2_L'+\eqref{new})$, respectively. Note that, again, we have included valid inequalities \eqref{desthub} and \eqref{orhub}. Besides, in the column related to the gap at termination, we report in parenthesis the number of instances (if any) that were not solved within the time limit.

Observe that these formulations allow to solve most of the considered instances in less than two hours. Observe that formulations $(F1_L)$+\eqref{flow1}-\eqref{flow3} and $(F1_L)$+\eqref{sec} provide the best time performances. Particularly, $(F1_L)$+\eqref{sec} is able to solve more instances than the remaining formulations for $p=7$, $\alpha=0.2$ and (\% edges)$=0.6$.

\begin{sidewaystable}[htbp]
  \centering
  \caption{Time results of formulations $(F1_L)$+\eqref{flow1}-\eqref{flow3}, $(F1_L)$+\eqref{sec}, $(F2_L)$,$(F2'_L)$ and $(F2'_L)+\eqref{new}$,  with valid inequalities \eqref{desthub}+\eqref{orhub} for $n=25$}
   \begin{adjustbox}{max width=1\textwidth}
    \begin{tabular}{cclrrrrrrrrrrrrrrrrr}
       \toprule
    \multirow{2}[2]{*}{$p$} & \multirow{2}[2]{*}{$\alpha$} & \multicolumn{1}{c}{\multirow{2}[2]{*}{(\% edges)}} & \multicolumn{1}{c}{\multirow{2}[2]{*}{$n_{path}$}} & \multicolumn{1}{c}{\multirow{2}[2]{*}{$t_{path}$}} & \multicolumn{3}{c}{$(F1_L)$+\eqref{flow1}-\eqref{flow3}} & \multicolumn{3}{c}{$(F1_L)$+\eqref{sec}} & \multicolumn{3}{c}{$(F2_L)$} & \multicolumn{3}{c}{$(F2'_L)$} & \multicolumn{3}{c}{$(F2_L)$+\eqref{new}} \\
      \cmidrule(lr){6-8}
      \cmidrule(lr){9-11}
      \cmidrule(lr){12-14}
      \cmidrule(lr){15-17}
      \cmidrule(lr){18-20}
          &       &       &       &       & \multicolumn{1}{c}{Gap$_{lp}$} & \multicolumn{1}{c}{Gap} & \multicolumn{1}{c}{$t$} & \multicolumn{1}{c}{Gap$_{lp}$} & \multicolumn{1}{c}{Gap} & \multicolumn{1}{c}{$t$} & \multicolumn{1}{c}{Gap$_{lp}$} & \multicolumn{1}{c}{Gap} & \multicolumn{1}{c}{$t$} & \multicolumn{1}{c}{Gap$_{lp}$} & \multicolumn{1}{c}{Gap} & \multicolumn{1}{c}{$t$} & \multicolumn{1}{c}{Gap$_{lp}$} & \multicolumn{1}{c}{Gap} & \multicolumn{1}{c}{$t$} \\
       \midrule
    \multirow{3}[2]{*}{3} & \multirow{3}[2]{*}{0.8} & 0.4   & 1683  & 2.79  & 3.75  & 0.00  & 1.63  & 5.16  & 0.00  & \textbf{1.52} & 5.16  & 0.00  & 1.88  & 5.16  & 0.00  & 1.95  & 3.70  & 0.00  & 1.92 \\
          &       & 0.5   & 2366  & 3.23  & 4.51  & 0.00  & 2.03  & 5.58  & 0.00  & \textbf{1.86} & 5.58  & 0.00  & 2.85  & 5.58  & 0.00  & 3.22  & 4.43  & 0.00  & 2.63 \\
          &       & 0.6   & 3232  & 3.37  & 2.98  & 0.00  & 2.16  & 3.63  & 0.00  & \textbf{2.09} & 3.63  & 0.00  & 2.88  & 3.63  & 0.00  & 2.94  & 2.86  & 0.00  & 2.53 \\
    \midrule
    \multirow{3}[2]{*}{3} & \multirow{3}[2]{*}{0.5} & 0.4   & 12890 & 5.64  & 0.15  & 0.00  & 4.80  & 0.15  & 0.00  & \textbf{4.52} & 0.15  & 0.00  & 5.13  & 0.15  & 0.00  & 5.31  & 0.15  & 0.00  & 5.48 \\
          &       & 0.5   & 18758 & 6.27  & 0.00  & 0.00  & 5.67  & 0.00  & 0.00  & \textbf{5.35} & 0.00  & 0.00  & 6.26  & 0.00  & 0.00  & 6.16  & 0.00  & 0.00  & 6.00 \\
          &       & 0.6   & 25486 & 6.86  & 0.00  & 0.00  & 6.65  & 0.00  & 0.00  & \textbf{6.37} & 0.00  & 0.00  & 7.00  & 0.00  & 0.00  & 7.35  & 0.00  & 0.00  & 7.42 \\
    \midrule
    \multirow{3}[2]{*}{3} & \multirow{3}[2]{*}{0.2} & 0.4   & 39075 & 8.38  & 0.00  & 0.00  & 9.00  & 0.00  & 0.00  & \textbf{8.75} & 0.00  & 0.00  & 9.69  & 0.00  & 0.00  & 12.93 & 0.00  & 0.00  & 9.84 \\
          &       & 0.5   & 59291 & 10.08 & 0.00  & 0.00  & 12.20 & 0.00  & 0.00  & \textbf{12.10} & 0.00  & 0.00  & 12.92 & 0.00  & 0.00  & 13.23 & 0.00  & 0.00  & 12.91 \\
          &       & 0.6   & 81462 & 11.94 & 0.00  & 0.00  & 16.13 & 0.00  & 0.00  & \textbf{15.46} & 0.00  & 0.00  & 16.61 & 0.00  & 0.00  & 17.00 & 0.00  & 0.00  & 16.74 \\
    \midrule
    \multirow{3}[2]{*}{5} & \multirow{3}[2]{*}{0.8} & 0.4   & 2738  & 21.01 & 5.66  & 0.00  & \textbf{2.19} & 5.69  & 0.00  & 4.31  & 5.69  & 0.00  & 9.34  & 5.69  & 0.00  & 8.49  & 5.37  & 0.00  & 8.80 \\
          &       & 0.5   & 4611  & 45.73 & 6.12  & 0.00  & \textbf{3.38} & 6.12  & 0.00  & 8.34  & 6.12  & 0.00  & 12.15 & 6.12  & 0.00  & 9.29  & 6.12  & 0.00  & 16.13 \\
          &       & 0.6   & 7368  & 90.08 & 4.80  & 0.00  & \textbf{5.02} & 4.80  & 0.00  & 5.93  & 4.80  & 0.00  & 5.71  & 4.80  & 0.00  & 6.29  & 4.80  & 0.00  & 9.02 \\
    \midrule
    \multirow{3}[2]{*}{5} & \multirow{3}[2]{*}{0.5} & 0.4   & 46765 & 36.19 & 3.67  & 0.00  & \textbf{31.34} & 3.67  & 0.00  & 41.64 & 3.67  & 0.00  & 112.38 & 3.67  & 0.00  & 78.76 & 3.67  & 0.00  & 116.83 \\
          &       & 0.5   & 103808 & 76.02 & 4.08  & 0.00  & \textbf{78.81} & 4.08  & 0.00  & 109.02 & 4.08  & 0.00  & 170.58 & 4.08  & 0.00  & 127.81 & 4.08  & 0.00  & 203.76 \\
          &       & 0.6   & 181026 & 151.58 & 3.40  & 0.00  & 121.64 & 3.40  & 0.00  & \textbf{111.35} & 3.40  & 0.00  & 217.15 & 3.40  & 0.00  & 356.19 & 3.40  & 0.00  & 430.03 \\
    \midrule
    \multirow{3}[2]{*}{5} & \multirow{3}[2]{*}{0.2} & 0.4   & 624037 & 131.03 & 3.03  & 0.00  & 615.54 & 3.03  & 0.00  & \textbf{420.70} & 3.03  & 0.01  & 2034.07 & 3.03  & 0.00  & 975.48 & 3.03  & 0.00  & 1470.03 \\
          &       & 0.5   & 1553479 & 308.46 & 0.42  & 0.00  & 545.70 & 0.42  & 0.00  & \textbf{477.05} & 0.42  & 0.00  & 560.85 & 0.42  & 0.00  & 1057.20 & 0.42  & 0.00  & 687.97 \\
          &       & 0.6   & 3019015 & 601.96 & 1.01  & 0.28  & (1)3001.82 & 1.01  & 0.00  & 2416.41 & 1.01  & 0.26  & (1)3010.79 & 1.01  & 0.00  & \textbf{2395.10} & 1.01  & 0.00  & 3287.66 \\
    \midrule
    \multirow{3}[2]{*}{7} & \multirow{3}[2]{*}{0.8} & 0.4   & 2772  & 957.35 & 5.18  & 0.00  & \textbf{3.01} & 5.24  & 0.00  & 10.05 & 5.24  & 0.00  & 22.85 & 5.24  & 0.00  & 19.82 & 5.19  & 0.00  & 14.00 \\
          &       & 0.5   & 4719  & 3654.74 & 8.28  & 0.00  & \textbf{11.49} & 8.28  & 0.00  & 19.66 & 8.28  & 0.00  & 50.41 & 8.28  & 0.00  & 30.03 & 8.28  & 0.00  & 39.91 \\
          &       & 0.6   & 7643  & 13657.80 & 9.02  & 0.00  & \textbf{25.50} & 9.02  & 0.00  & 27.16 & 9.02  & 0.00  & 60.65 & 9.02  & 0.00  & 45.78 & 9.02  & 0.00  & 52.51 \\
    \midrule
    \multirow{3}[2]{*}{7} & \multirow{3}[2]{*}{0.5} & 0.4   & 82322 & 1405.74 & 4.77  & 0.00  & 136.64 & 4.77  & 0.00  & 276.80 & 4.77  & 0.00  & 379.34 & 4.77  & 0.00  & \textbf{114.88} & 4.77  & 0.00  & 340.60 \\
          &       & 0.5   & 250570 & 5373.92 & 6.98  & 0.00  & 600.35 & 6.98  & 0.00  & 808.67 & 6.98  & 0.01  & 1730.97 & 6.98  & 0.01  & \textbf{379.59} & 6.98  & 0.01  & 1428.58 \\
          &       & 0.6   & 552628 & 16295.21 & 7.08  & 0.00  & \textbf{887.65} & 7.08  & 0.03  & (1)2700.30 & 7.60  & 1.89  & (2)5153.91 & 7.08  & 0.17  & (1)2567.73 & 7.08  & 0.86  & (1)5062.37 \\
    \midrule
    \multirow{2}[2]{*}{7} & \multirow{2}[2]{*}{0.2} & 0.4   & 8198377 & 4163.93 & 7.42  & 5.56  & (4)8532.19 & 7.48  & 5.19  & \textbf{(4)8634.12} & 7.61  & 6.88  & (5)8657.80 & 7.56  & 6.67  & (5)7887.83 & 76.20 & 75.50 & (5)8663.45 \\
          &       & 0.5   & 33513751 & 16063.98 & $>$100  & $>$100  & (5)13573.25 & 24.05 & 24.05 & \textbf{(5)12901.71} & $>$100  & $>$100  & (5)13003.63 & $>$100  & $>$100  & (5)10663.05 & $>$100  & $>$100  & (5)13290.84 \\
    \bottomrule
    \end{tabular}%
    \end{adjustbox}
  \label{tab_25}%
\end{sidewaystable}%

Table \ref{tab_25} follows a similar structure to Table \ref{tab_20}. We use $n=25$, $p\in\{3,5,7\}$, $\alpha\in\{0.2,0.5,0.8\}$ and (\%edges)$\in\{0.4,0.5,0.6\}$. 
Observe that, since we are considering a maximization problem, in some cases, we obtain gap percentages greater than $100\%$. In those cases, we report ``$>$100'' in the table. Besides, for $n=25$, $p=7$, $\alpha=0.2$ and $(\%edges)=0.6$ we obtained out-of-memory status, consequently,  Table \ref{tab_25} does not include the results for these parameter values.

Observe that, in Table \ref{tab_25}, the largest average number of paths is 33,513,751 that are obtained in an average of 16,063.98 seconds. In addition, formulation $(F1_{L})$+\eqref{sec} provides the best time results.  Note that, in general, this formulation can solve most of the instances within the time limit, and it is the one providing the best gaps for $p=7$ and $\alpha=0.2$.

Overall, the results of tables \ref{tab_nonlinear} to \ref{tab_25} show that MILP formulations outperform the results of the MINLP formulations. Although it is necessary to compute the set of candidate paths to use MILP formulations, the efficient algorithm provided in Section \ref{paths} allows us to solve instances with $n=10$, $n=15$, $n=20$ and $n=25$.

\subsection{Study case: determining a hub line in Montreal}\label{subsecMontreal} 

In this section, we select Montreal (Canada) as our study case to illustrate the application of the proposed model. We analyze the impact of the parameter  $r$ considered in the gravity model and the percentage $\vartheta$ used for the access and exit times. Then, we compare the configuration of the hub network using elastic (ED-HLLP)  with respect to the network using  inelastic demand (HLLP).

\subsubsection{Data description} \label{data-montreal}

Montreal is the second-largest city in Canada, and thousands of tourists travel to this city every day. Therefore, designing an efficient public transportation system for this population is of utmost importance. Montreal is divided into 40 different boroughs based on data available from Statistics Canada, known as census metropolitan areas (CMA). The centroids of these boroughs will be the considered  origins and destinations nodes in our problem. 

\begin{figure}[h]
  \centering
  \includegraphics[scale=0.35]{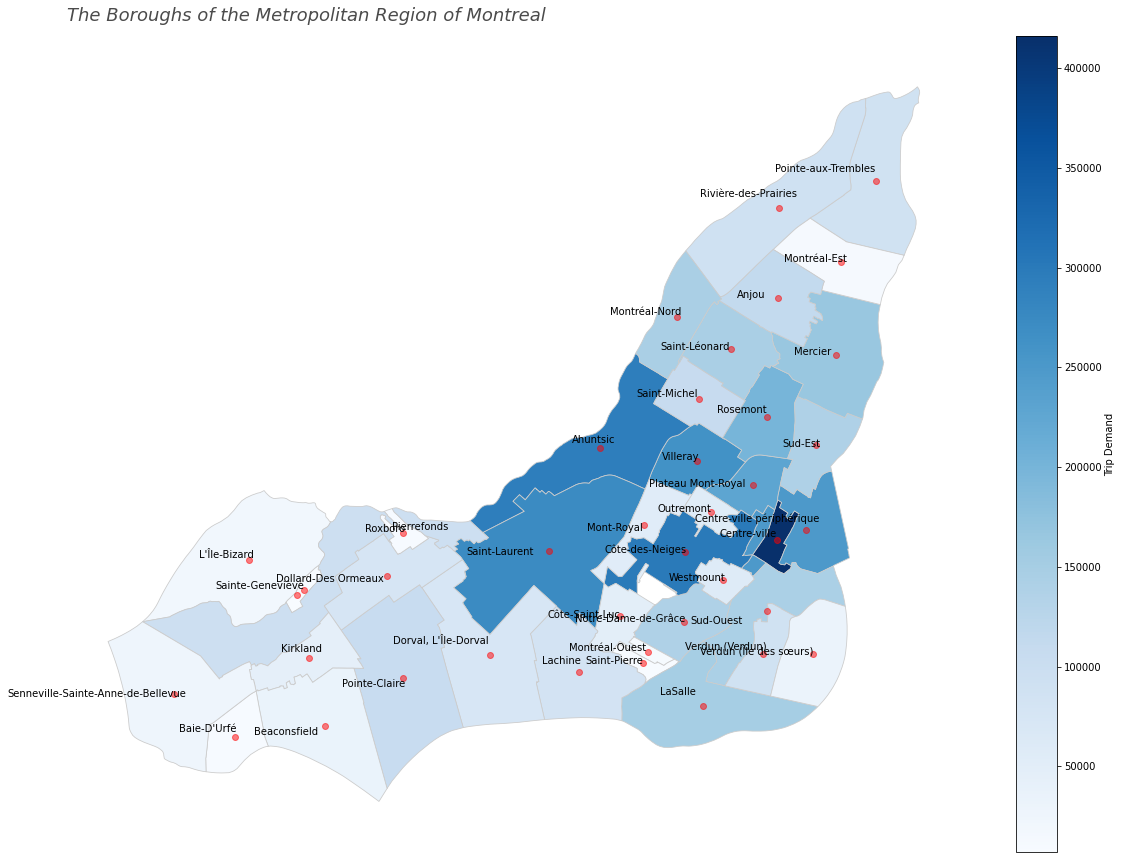}
   \caption{ The boroughs of the Metropolitan Region of Montreal}
  \label{montreal}
\end{figure}

Figure \ref{montreal} shows the locations and the corresponding name of the nodes on a map of Montreal; the map was plotted using  the GeoPandas library from Python 3.8. Besides, the terrestrial limits of each CMA of the Montreal were obtained from 2016 census available in \url{https://www12.statcan.gc.ca/census-recensement/2011/geo/bound-limit/bound-limit-2016-eng.cfm}.

For this analysis, the 2018 OD survey of Montreal was used. This survey has been conducted every five years since 1970. The Autorit\'e R\'egionale de Transport M\'etropolitain (ARTM) produces this OD survey in partnership with the Minister of Transports of Qu\'ebec, the R\'eseau de transport de Longueuil, the Society of transport de Montréal (STM), the Society de transport de Laval, and the Montreal Metropolitan Community.
It is a complete source of information on the movement of people on foot, bike, bus, metro, train or car, in the metropolitan area of Montreal. For more details, we refer to \cite{web:ARTM}. 

The times between boroughs' centroids were obtained using Google Maps API in Python 3.8. In addition, a post-processing was carried out in the times obtained between each node in order to ensure the triangle inequality is satisfied. The instances used in this study can be found in the following repository, \textit{https://github.com/brenda} \textit{denisse16/HLLP-ELASTIC-DEMAND/tree/main/data.}

\subsubsection{Sensitivity analysis of $r$ and $\vartheta$} \label{experiments-montreal}

In this section, two different sensitivity analyses are carried out about how the values of parameters $r$ and $\vartheta$ affect the establishment of the hub network in Montreal. For these studies we considered different numbers of nodes ($n$). Particularly, the nodes related to the neighborhoods with the highest travel demands were selected. Furthermore, different values of $p$, $r$ and $\vartheta$ were used.

\begin{footnotesize}
\begin{table}[h!]
\centering
\caption{Different values of $r$}
\label{r_values}
\begin{tabular}{@{}lr@{}}
\toprule
transport mode (reference)            & \multicolumn{1}{c}{$r$} \\ \midrule
Railway express \cite[][]{zipf1946p}  & 1.00          \\ \midrule
Train flows \cite[][]{sun2019spatial} &  1.70 \\
\midrule
Subway \cite[][]{goh2012modification}         & 2.68       \\
\bottomrule
\end{tabular}
\end{table}
\end{footnotesize}

\begin{table}[h!]
\centering
\caption{Obtained results when changing $r$ value}
\label{r_montreal}
\begin{tabular}{ccccccc}
\hline
$n$ &
  $p$ &
 $\vartheta$ &
  $r$ &
\begin{tabular}[c]{@{}c@{}}\% OD pairs  \\  served\end{tabular} &
  \begin{tabular}[c]{@{}c@{}}\% served \\ demand\end{tabular} &
  \begin{tabular}[c]{@{}c@{}}\% saved \\ time\end{tabular} \\ \hline
20 & 5 & 0.1 & 1.00 & 12 & 29 & 31 \\ \hline
20 & 5 & 0.1 & 1.70 & 14 & 37 & 27 \\ \hline
20 & 5 & 0.1 & 2.68 & 15 & 51 & 21 \\ \hline
\end{tabular}
\end{table}


\begin{figure}[h!]
	\centering
	\begin{subfigure}[b]{0.55\textwidth}
		\centering
		\includegraphics[width=1.1\linewidth]{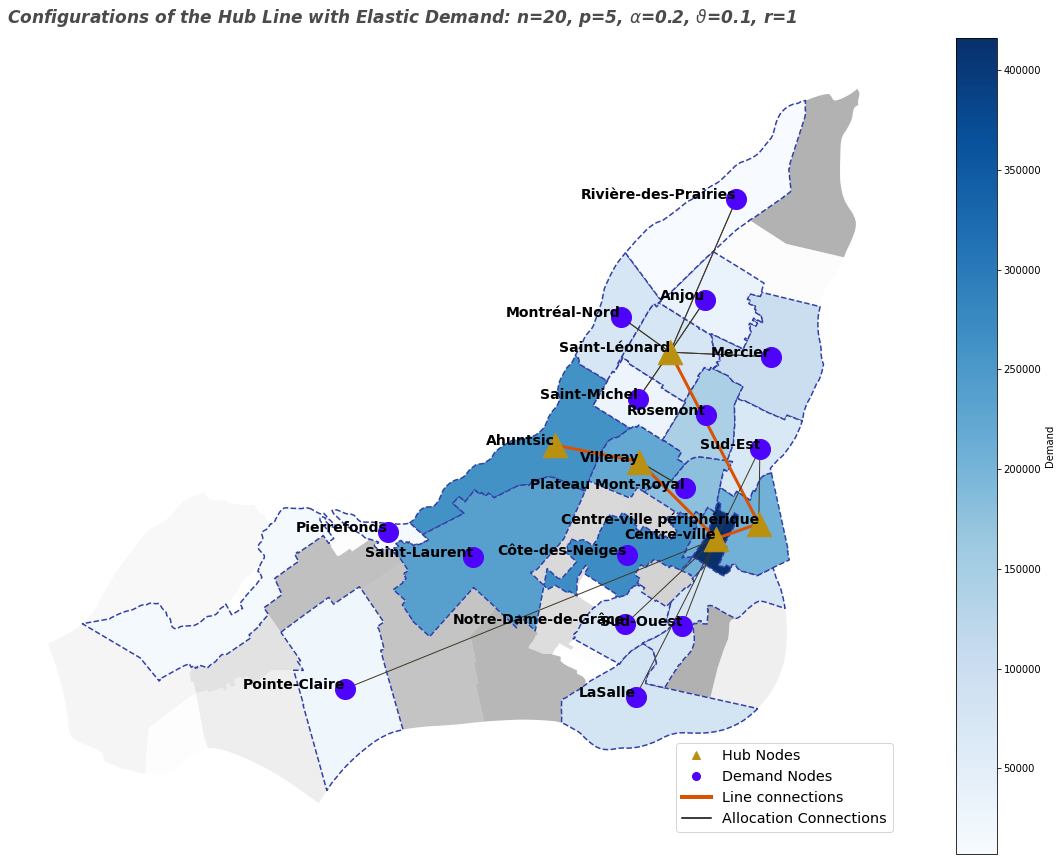}
		\caption{$r$=1 }
		\label{hubn20r1}
	\end{subfigure}%
	\begin{subfigure}[b]{0.55\textwidth}
		\centering
		\includegraphics[width=1.1\linewidth]{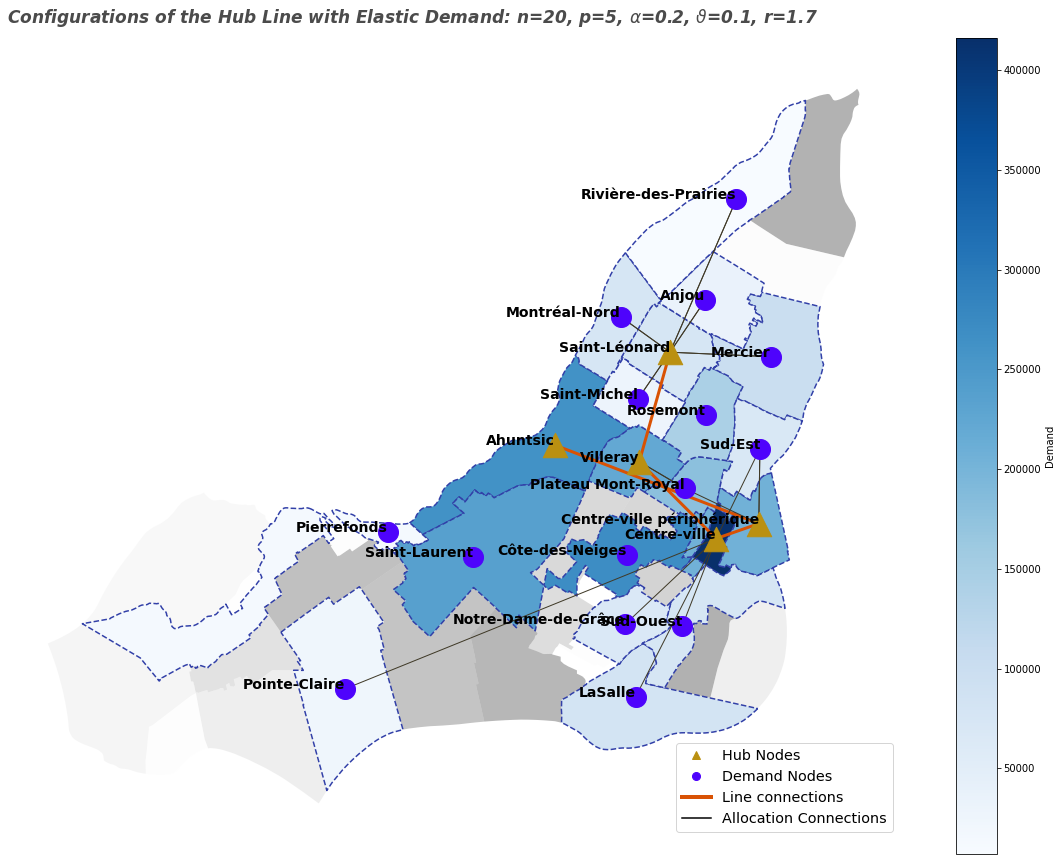}
		\caption{ $r$=1.7}
		\label{hubn20r17}
	\end{subfigure}%
	\newline
	\begin{subfigure}[b]{0.55\textwidth}
		\centering
		\includegraphics[width=1.1\linewidth]{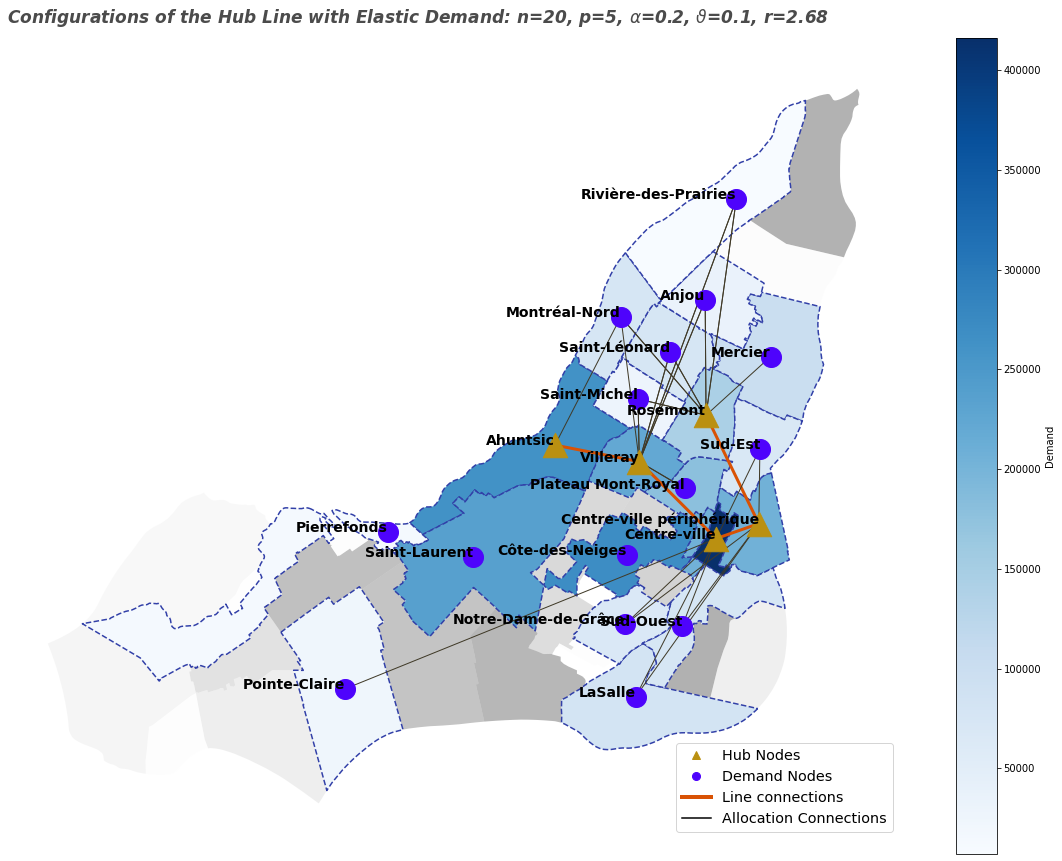}
		\caption{$r$=2.68}
		\label{hubn20r268}
	\end{subfigure}
	\caption{Configuration of the hub line 
	for $n$=20, $p$=5, $\alpha$ =0.2,$\vartheta$ =0.1 and different values of $r$.}
\end{figure}

Recall that parameter $r$ is the exponent considered in the gravity model, and it is related to the efficiency of the selected transport mode. For this analysis, we considered three different values of $r$. Particularly, we selected $r$ values that have been previously studied in the literature for different transport modes.  In Table \ref{r_values}, we report the considered values of $r$ and some references about the use of these specific values.
    
Table \ref{r_montreal} and Figures \ref{hubn20r1}, \ref{hubn20r17} and \ref{hubn20r268} show the results provided by the model when using these values of $r$. 
Table \ref{r_montreal} reports the values of parameters used for the model in the first four columns, the fifth column shows the percentage of OD pairs that use the hub line (\% OD pairs served), the sixth column presents the percentage of total demand which is transported using the hub line (\% served demand), and finally, the last column shows the average saved time after the establishment of the hub line (\% saved time). The percentage of saved time is the rate of the difference between the total travel time before and after the location of the hub line divided by the former.

Observe that, when $r$ increases, the percentage of ODs using the hub line and the percentage of the demand transported through the hub line also increase. In Figure \ref{hubn20r268}, when $r$ = 2.68, in the surrounding of the selected hubs there are big companies, shopping centers, etc. For example, {\sl Centre-Ville} and {\sl Centre-ville-P\'eriph\'erique} hubs are located in the downtown Montreal, and they include the central business district and main commercial streets. Also, $Rosemont$ borough has important points of interest (POI) such as Jean-Talon Market, Maisonneuve Park, Montreal Botanical Garden and Montreal Public Libraries. 

Although the percentage of saved time decreases as the value of $r$ increases, we can see that the selected hubs are located in areas with a high concentration of urban activity and a high demand flow.

\begin{table}[h!]
\centering
\caption{Results of changing the $\vartheta$ }
\label{montreal_v}
\begin{tabular}{@{}cccccc@{}}
\toprule
\rowcolor[HTML]{FFFFFF} 
$n$ &
$p$ &
 $\vartheta$ &
  \begin{tabular}[c]{@{}c@{}}\% O/D pairs \\ served\end{tabular} &
  \begin{tabular}[c]{@{}c@{}}\% served\\  demand\end{tabular} &
  \begin{tabular}[c]{@{}c@{}}\% saved\\  time \end{tabular} \\ \midrule
\textbf{15} & \textbf{7} & 0.00 & 51 & 98 & 33\\ \midrule
\textbf{15} & \textbf{7} & 0.10 & 29 & 68 & 27 \\ \midrule
\textbf{15} & \textbf{7} & 0.15 & 21 & 53 & 26 \\ \midrule
\textbf{15} & \textbf{7} & 0.25 & 16 & 20 & 20 \\ \bottomrule
\end{tabular}
\end{table}

\begin{figure}[h!]
	\centering
	\begin{subfigure}[b]{0.55\textwidth}
		\centering
		\includegraphics[width=1.0\linewidth]{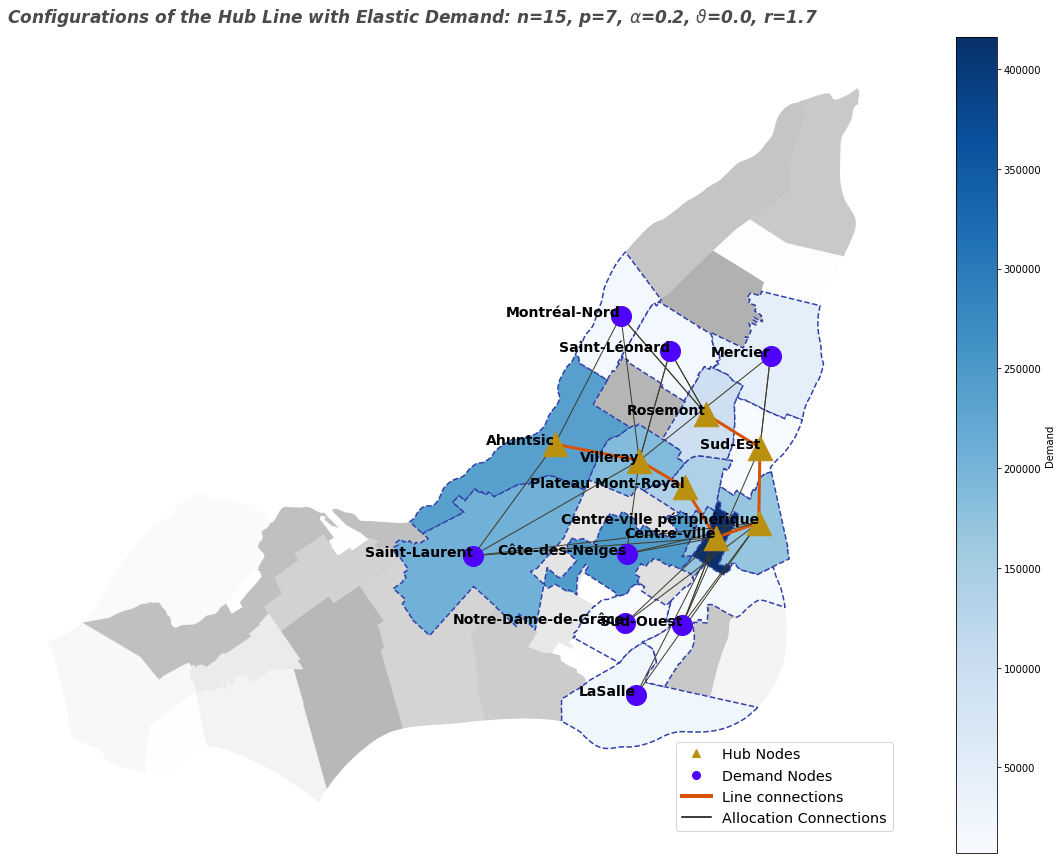}
		\caption{$\vartheta$=0 }
		\label{hubn15v0}
	\end{subfigure}%
	\begin{subfigure}[b]{0.55\textwidth}
		\centering
		\includegraphics[width=1.0\linewidth]{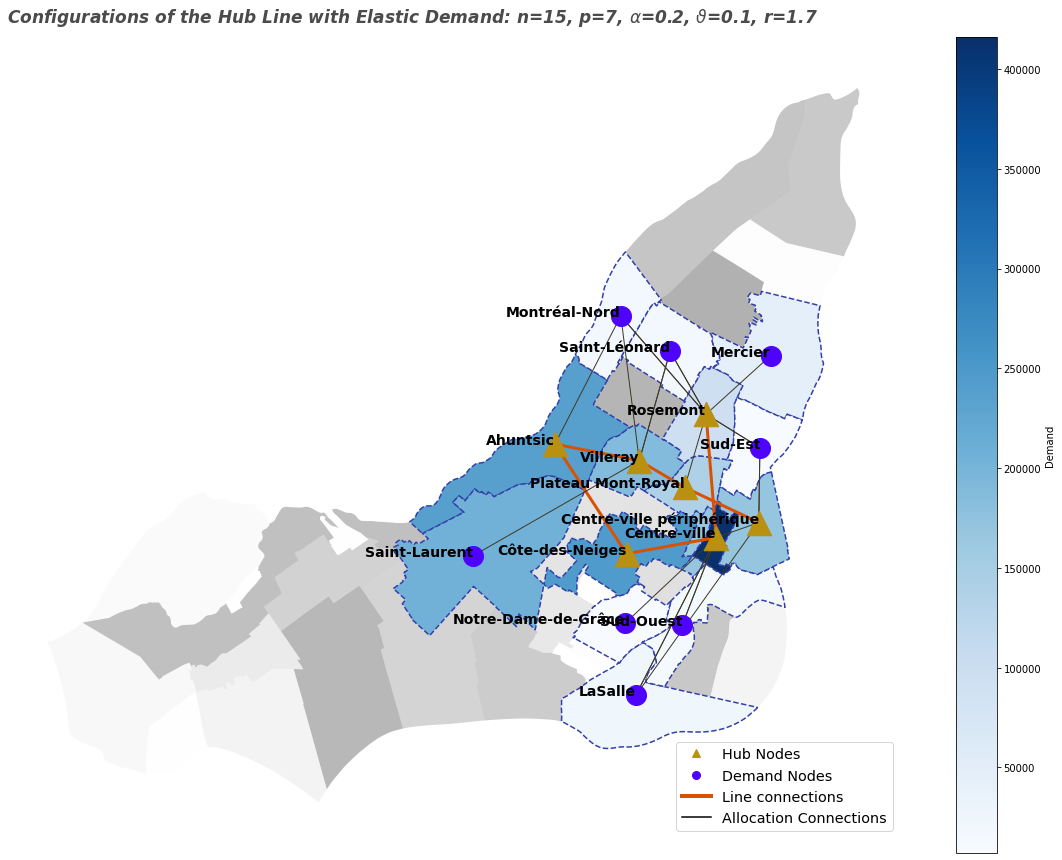}
		\caption{$\vartheta$=0.1}
		\label{hubn15v01}
	\end{subfigure}%

	\begin{subfigure}[b]{0.55\textwidth}
		\centering
		\includegraphics[width=1.0\linewidth]{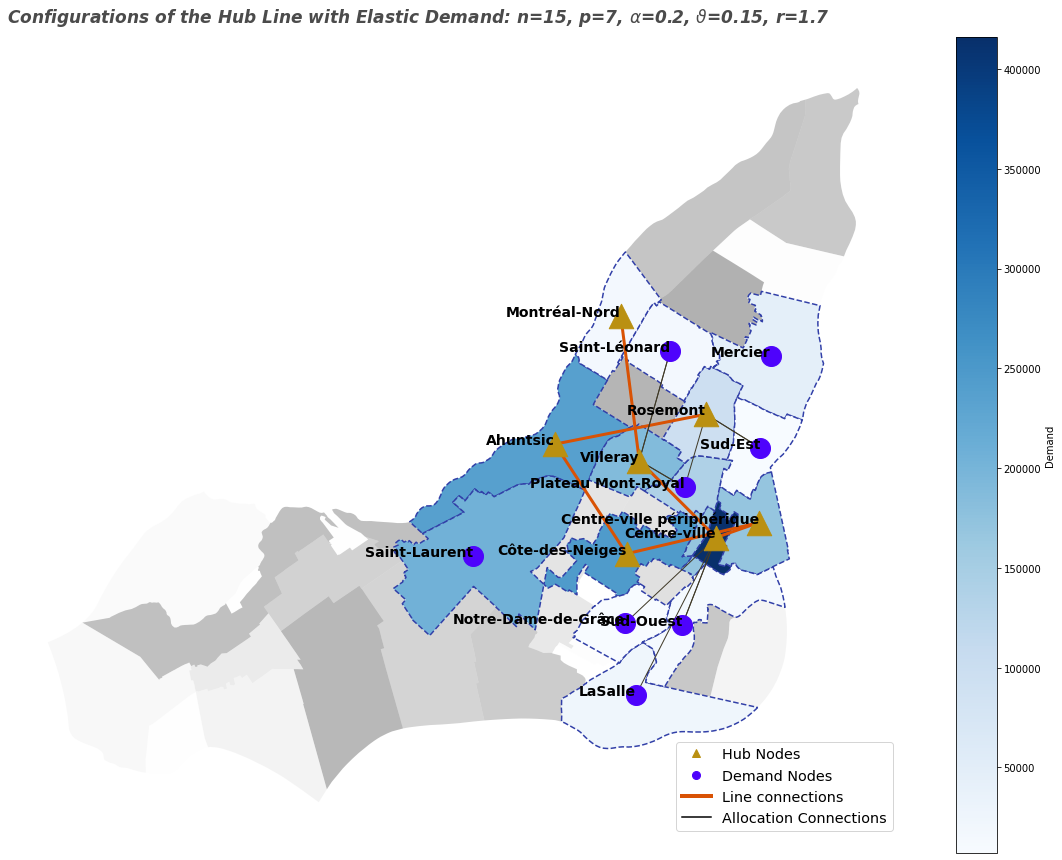}
		\caption{$\vartheta$=0.15}
		\label{hubn15v015}
	\end{subfigure}%
	\begin{subfigure}[b]{0.55\textwidth}
	\centering
	\includegraphics[width=1.0\linewidth]{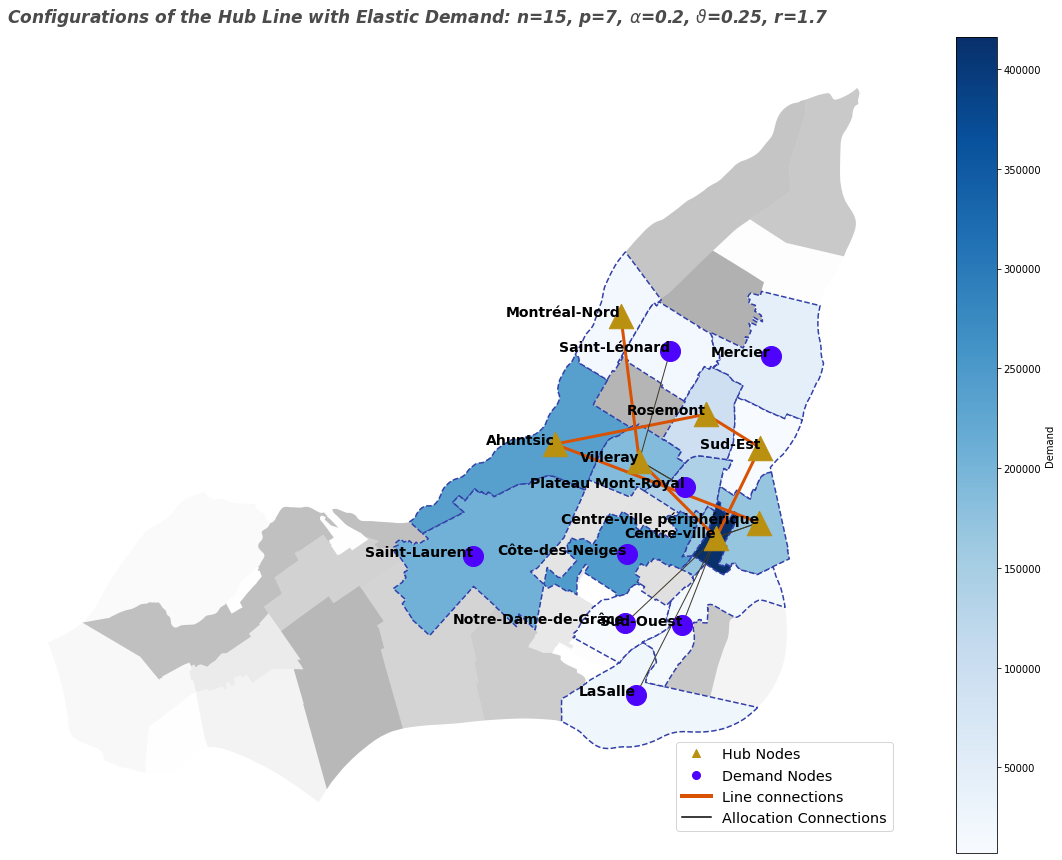}
	\caption{$\vartheta$=0.25}
	\label{hubn15v025}
	\end{subfigure}%
	
	\caption{Configuration of the hub line for $n$=15, $p$=7, $\alpha$=0.2, $r$=1.7, and different values of $\vartheta$.}
\end{figure}

We also analyze the performance of the \ED for different values of $\vartheta$. We consider $\vartheta\in \{0, 0.1, 0.15, 0.25\}$. Recall that this parameter is used to determine the access and exit times. For instance, setting $\vartheta$=0.15, means that the access and exit time are 15\% of the average travel time. 

Table \ref{montreal_v} and Figures from \ref{hubn15v0} to \ref{hubn15v025}  show the different configurations obtained by varying $\vartheta$. Table \ref{montreal_v} follows the same configuration as Table \ref{r_montreal}. Note that the percentage of OD pairs using the hub line, the percentage of total demand transported using the line, and the average saved time decrease when $\vartheta$ increases. Consequently, we appreciate that big access and exit times make the use of hub network less appealing. 

\subsubsection{Comparison of \ED with respect to the determination of the hub line considering inelastic demand} \label{elastic_noelastic_montreal}

\begin{table}[h!]
\centering
\caption{Comparison of ED-HLLP and HLLP}
\label{comparison}
\begin{adjustbox}{max width=1\textwidth}
\begin{tabular}{cccrrrrrr}
\hline
\multirow{2}{*}{$n$} &
  \multirow{2}{*}{$p$} &
  \multirow{2}{*}{ $\vartheta$} &
  \multicolumn{2}{c}{\begin{tabular}[c]{@{}c@{}}\%  O/D pairs \\ served\end{tabular}} &
  \multicolumn{2}{c}{\begin{tabular}[c]{@{}c@{}}\% served \\ demand\end{tabular}} &
  \multicolumn{2}{c}{\begin{tabular}[c]{@{}c@{}}\% saved time\\      \end{tabular}} \\ 
    \cmidrule(lr){4-5}
    \cmidrule(lr){6-7}
    \cmidrule(lr){8-9}
 &
   &
   &
  \multicolumn{1}{c}{Elastic} &
  Static &
  \multicolumn{1}{c}{Elastic} &
  Static &
  \multicolumn{1}{c}{Elastic} &
  Static \\ 
\midrule 
  \multirow{3}[2]{*}{10} & \multirow{3}[2]{*}{3} & 0.10  & \multicolumn{1}{r}{16}          & 16 & \multicolumn{1}{r}{36} & 19 & \multicolumn{1}{r}{27}          & 31 \\ 
 &  & 0.15 & \multicolumn{1}{r}{11}          & 13 & \multicolumn{1}{r}{25} & 17 & \multicolumn{1}{r}{26}          & 26 \\ 
 &  & 0.25 & \multicolumn{1}{r}{7}           & 8  & \multicolumn{1}{r}{14} & 8  & \multicolumn{1}{r}{21}          & 24 \\ 
\midrule
  \multirow{3}[2]{*}{10} & \multirow{3}[2]{*}{5}  & 0.10  & \multicolumn{1}{r}{29} & 28 & \multicolumn{1}{r}{64} & 43 & \multicolumn{1}{r}{37}          & 40 \\ 
 &  & 0.15 & \multicolumn{1}{r}{27}          & 28 & \multicolumn{1}{r}{53} & 31 & \multicolumn{1}{r}{30}          & 30 \\ 
 &  & 0.25 & \multicolumn{1}{r}{11}          & 16 & \multicolumn{1}{r}{21} & 20 & \multicolumn{1}{r}{28} & 26 \\ \midrule
  \multirow{3}[2]{*}{10} & \multirow{3}[2]{*}{7}  & 0.10  & \multicolumn{1}{r}{43}          & 52 & \multicolumn{1}{r}{75} & 60 & \multicolumn{1}{r}{37}          & 37 \\ 
 &  & 0.15 & \multicolumn{1}{r}{32}          & 46 & \multicolumn{1}{r}{56} & 59 & \multicolumn{1}{r}{32} & 31 \\ 
 &  & 0.25 & \multicolumn{1}{r}{24}          & 30 & \multicolumn{1}{r}{36} & 36 & \multicolumn{1}{r}{24}          & 21 \\ \midrule
  \multirow{3}[2]{*}{15} & \multirow{3}[2]{*}{3}  & 0.10 & \multicolumn{1}{r}{12} & 11 & \multicolumn{1}{r}{22} & 14 & \multicolumn{1}{r}{21}          & 23 \\ 
 &  & 0.15 & \multicolumn{1}{r}{9}           & 10 & \multicolumn{1}{r}{15} & 7  & \multicolumn{1}{r}{20}          & 22 \\ 
 &  & 0.25 & \multicolumn{1}{r}{6}           & 5  & \multicolumn{1}{r}{6}  & 5  & \multicolumn{1}{r}{19}          & 23 \\ \midrule
  \multirow{3}[2]{*}{15} & \multirow{3}[2]{*}{5}  & 0.10  & \multicolumn{1}{r}{21}          & 23 & \multicolumn{1}{r}{48} & 23 & \multicolumn{1}{r}{26}          & 29 \\ 
 &  & 0.15 & \multicolumn{1}{r}{15}          & 17 & \multicolumn{1}{r}{30} & 21 & \multicolumn{1}{r}{29} & 27 \\ 
 &  & 0.25 & \multicolumn{1}{r}{10}          & 10 & \multicolumn{1}{r}{13} & 12 & \multicolumn{1}{r}{23}          & 23 \\ \midrule
  \multirow{3}[2]{*}{15} & \multirow{3}[2]{*}{7}  & 0.10  & \multicolumn{1}{r}{29}          & 32 & \multicolumn{1}{r}{57} & 41 & \multicolumn{1}{r}{27}          & 29 \\ 
 &  & 0.15 & \multicolumn{1}{r}{23}          & 28 & \multicolumn{1}{r}{45} & 35 & 30 & 27 \\ 
 &  & 0.25 & 16 & 14 & \multicolumn{1}{r}{21} & 21 & \multicolumn{1}{r}{20}          & 21 \\ 
\bottomrule
\end{tabular}%
\end{adjustbox}
\end{table}

\begin{figure}[h!]
	\centering
	\begin{subfigure}[b]{0.55\textwidth}
		\centering
		\includegraphics[width=0.9\linewidth]{hub_line_elastic_F3DV_n15_p7_alpha0.2_v0.15_r1.7}
		\caption{elastic}
		\label{hubn15v015elastic}
	\end{subfigure}%
	\begin{subfigure}[b]{0.55\textwidth}
		\centering
		\includegraphics[width=0.9\linewidth]{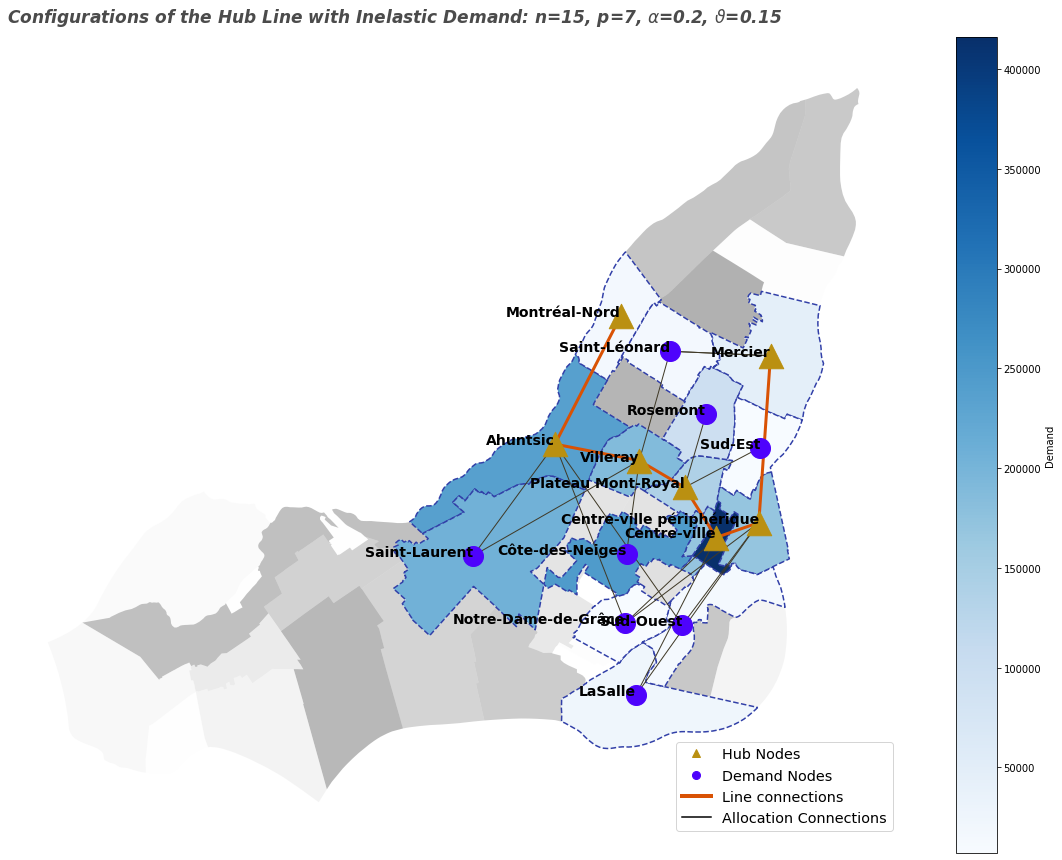}
		\caption{static}
		\label{hubn15v015inelastic}
	\end{subfigure}%
	
	\caption{Hub line configuration for the study case with elastic and static demand}
\end{figure}

Table \ref{comparison} shows the percentage of OD pairs using the hub line, the percentage of demand which is transported through the hub line and the percentage of saved time for the elastic and inelastic version of the hub line location problem. The results of the case study with $n\in\{10,15\}$, $p\in\{3,5,7\}$, $\vartheta \in \{0.1, 0.15, 0.25\}$, $\alpha=0.2$ and $r$=1.7 are reported in this table.

Regarding the percentage of commodities using the hub line, we can observe, that in general, more commodities use the line when considering inelastic demand. However, in some cases better results are obtained with the proposed model. Besides, we can conclude that both models provide similar percentages of saved times. 

Observe that the percentage of demand which is transported through the hub line is bigger when applying the \ED proposed in this paper. For instance, for $n=15$, $p=7$ and $\vartheta=0.7$, the 57\% of demand is served using the hub line, while considering static demand only 41\% of the demand is routed through the line.

Figures \ref{hubn15v015elastic} and \ref{hubn15v015inelastic} show the obtained hub lines when using elastic and inelastic demand, respectively. Note that proposed lines are very different. In particular, we can see in the in Figure \ref{hubn15v015elastic} that one of the hubs selected is ``C\^ote-des-Neiges'', which is a borough that has important POIs such as Montreal University, Hospitals and Saint Joseph's Oratory.

\section{Conclusions}
\label{conclusions}

The hub line location problem introduced in \cite{DeSa2015} aims to locate $p$ hubs connected in a line composed by $p-1$ hub edges in such a way that the total travel time associated with the commodities is minimized. In this paper, we proposed and analyzed an extension of this problem: the profit oriented hub line location problem with elastic demands (\ED). This new model introduces an objective function based on gravity models to include demand elasticity. The objective of this model is to maximize the total reduction time of the commodities when using the hub line.

We have proposed two main MINLP formulations ($F1$ and $F2$) to deal with this problem. Besides, we have developed a preprocessing phase, some fixing variables criteria and some valid inequalities. The properties of the nonlinear objective function have also been described.

In addition, we have introduced two main MILP formulations ($F1_L$ and $F2_L$) to address the problem. These MILP formulations assume that all the candidate paths for the commodities are known. In order to obtain all path candidates, we have developed an algorithm which takes into account the time reductions appearing in hub lines.

The computational results showed that formulations $(F1_L)$, $(F2_L)$, $(F2_L')$ and $(F2_L')$+\eqref{new} outperform the running times of ($F1$), ($F2$), $(F2')$ and ($F2'$+\eqref{new}). Besides, we can remark that the best time results are the ones corresponding to formulation $(F1_L)$ with constraints \eqref{sec} and valid inequalities \eqref{desthub} and \eqref{orhub}. Note also that Algorithm \ref{alg1} provides all candidate paths in suitable running times. 

 In addition, the presented model (\ED) was tested using real data from Montreal. Some comparisons have been reported considering different values of parameters $r$ and $\vartheta$. Also, we compared the results of \ED with respect to the hub line location problem with inelastic demand. 
 
 {Although this paper focuses on the design of a single hub line topology, additional constraints can be incorporated to avoid arc-crossings or to include preferences to reduce the travel times between extreme hub nodes (e.g., end-of-line terminals). We are also interested in further extending this work to consider multiple hub lines.}

\section*{Acknowledgments}

We thank Alex Guindon, GIS and Data Librarian/GPE Librarian at Concordia University, for his help in providing the sources of information for the Montreal data used in the study. Luisa I. Martínez-Merino and Antonio M. Rodríguez-Chía thank the Agencia Estatal de Investigación (AEI) and the European Regional Development’s funds (ERDF): projects TED2021-130875B-I00 and PID2020-114594GB-C22; Agencia Estatal de Investigación
(AEI): project MTM-RED2018-102363-T; Regional Government
of Andalusia: projects FEDER-UCA18-106895 and P18-FR-1422; and Fundación BBVA: project NetmeetData
(Ayudas Fundación BBVA a equipos de investigación científica 2019). Luisa I. Martínez-Merino acknowledges NSERC and Concordia (CARA-CASA OVPRGS), Subprograma Juan de la Cierva Formación 2019 (FJC2019-039023-I), PAIDI 2020 postdoctoral fellowship funded by
European Social Fund and Junta de Andalucía, project CIGE/2021/161 OPTICODS from Generalitat Valenciana, and poject AT21\_0003 2 from Junta de Andaluc\'ia. Ivan Contreras and Brenda Cobeña gratefully acknowledges the support by grant 2018–
06704 of the Canadian Natural Sciences and Engineering Research Council (NSERC). {We would like to thank the reviewers for their constructive suggestions to improve the
manuscript.}

\newpage


\end{document}